\documentclass[12pt]{article}
\usepackage{amsthm}
\usepackage{amsmath}
\usepackage{enumerate}
\usepackage[usenames,dvipsnames]{color}

  \newcommand{\lab}[1]{\label{#1}}
 \def\fullpage
{
\addtolength{\topmargin}{-1.5 cm}
\addtolength{\oddsidemargin}{-1.5 cm}
\addtolength{\textwidth}{+3 cm}
\addtolength{\textheight}{+3 cm}
}
\fullpage

\newtheorem{theorem}{Theorem}[section]
\newtheorem{lemma}[theorem]{Lemma}
\newtheorem{claim}[theorem]{Claim}
\newtheorem{corollary}[theorem]{Corollary}
\newtheorem{problem}[theorem]{Problem}

\newcommand{\bel}[1]{\be\lab{#1}}
\newcommand\eqn[1]{(\ref{#1})}
\newcommand{\be}{\begin{equation}}
\newcommand{\ee}{\end{equation}}
\newcommand{\bea}{\begin{eqnarray}}
\newcommand{\eea}{\end{eqnarray}}
\newcommand{\bean}{\begin{eqnarray*}}
\newcommand{\eean}{\end{eqnarray*}}
\def\blackslug{\hbox{\kern1pt\vrule height6pt width4pt depth1pt\kern1pt}}
\def\proof{\par\noindent{\bf Proof \enspace}\rm}
\def\qed{\penalty 500\hbox{\quad\blackslug}\ifmmode\else\par
                   \vskip4.5pt plus3pt minus2pt\fi}
\catcode`@=11 \@addtoreset{equation}{section}

\catcode`@=12

\def\aas{a.a.s.}

\def\B{{\cal B}}
\def\P{{\cal P}}
\def\G{{\cal G}}
\def\pr{{\bf P}}
\def\ex{{\bf E}}
\def\ctbar{\overline C_t}
\def\eps{{\epsilon}}
\def\la{{\lambda}}
\def\Bf{{B_{\rm final}}}
\def\Bfk{{B^k_{\rm final}}}
\def\Gf{{G_{\rm final}}}
\def\D{\overleftarrow{D_i}}
\def\Ds{\overleftarrow{D}}
\def\DD{D_i^\bullet}
\newcommand{\Db}{D^\bullet}
\title{Sharper analysis of the random graph $d$-process via a balls-in-bins model
\thanks{Most of this work was performed during mutual visits of the authors at, resp.,  Adam
Mickiewicz University (2009 and 2022) and Monash University (2014 and 2023).}}
\author{Andrzej Ruci{\' n}ski\thanks{Research supported in part by Narodowe Centrum Nauki, grant  2018/29/B/ST1/00426.
}\\ {\small Department of Discrete Mathematics}\\
{\small
Adam Mickiewicz University}\\
{\small Pozna\'n, Poland}\\
{\small {\tt rucinski@amu.edu.pl}} \and Nick Wormald\thanks{Research supported by the Canada Research Chairs program, NSERC, and the Australian Laureate Fellowships grant FL120100125. }
\\ {\small School of Mathematics}\\
 {\small Monash University}\\
 {\small VIC 3800 Australia}\\
  {\small {\tt nick.wormald@monash.edu }} }

\date{}

\begin{document}
\maketitle

\begin{abstract}
A graph $d$-process  starts with an empty graph on $n$ vertices, and adds one edge at each time
step, chosen uniformly at random from those pairs which are not yet edges and whose both vertices
have current degree less than $d$. If, in the final graph, at most one vertex has degree $d-1$ and all other have degree $d$, we call the process saturated.
We present a new approach to analysing this process based on
random allocation of balls in bins. This allows us to get improved results on the degree distribution throughout the process and, consequently, to determine the asymptotic probability of non-saturation of the process.

\end{abstract}

\noindent
\section{Introduction}\lab{s:intro}
Let $d$ and $n$ be positive integers. Consider a random process  $\P(d,n)$
 which adds edges, one
by one, to an originally empty graph on $n$ vertices, with the restriction that the next edge is
uniformly chosen out of all unused pairs containing vertices which have current degree less than
$d$ (called \emph{unsaturated}). Let us call the process the \emph{$d$-process}. After at most
$\lfloor  dn /2\rfloor$ steps, this process gets stuck at a final graph $\Gf$, which, however, does
not need to be $d$-regular, even when  $ dn $  is even.

We  say that an event holds \emph{asymptotically almost surely} (a.a.s.) if the probability that it
holds tends to 1 as $n\to\infty$. A $d$-process is said to \emph{saturate} if   $\Gf$
is either  $d$-regular or has $n-1$ vertices of degree $d$ and the remaining  vertex has degree
$d-1$.  Note that the $d$-process can only fail to saturate if the unsaturated vertices form a clique, and thus there are at most $d$ unsaturated vertices in $\Gf$.  In the 1980's, 'Erd{\H o}s  posed the problem (see~\cite{Ruc90}) of determining what  the limiting distribution of the number of unsaturated vertices is in $\Gf$. Answering that question, Ruci\'nski and Wormald~\cite{RW92} proved,  by a quite complex argument using martingales, that  the $d$-process  a.a.s.\
saturates.

\begin{theorem}[\cite{RW92}]\label{satur} For fixed $d\ge 2$, the $d$-process a.a.s.\ saturates.\end{theorem}

In this paper we propose a new approach to studying $\P(d,n)$ by linking it to a
balls-in-bins process which is much more amenable to study by an elementary approach. In this way we obtain arguably simpler  proofs of known
results, as well as new, more refined and extended results. In particular, some of our easier results (Lemma~\ref{l:X} and  Corollary~\ref{indept}) replace similar parts of the argument given in~\cite{RW92} using concentration via a Doob martingale, and enable  us to give a simpler proof of Theorem \ref{satur}. We are also able to prove the following refinement, which gives a more precise answer to Erd{\H o}s' original question.

 Let $F=F(n,d)$ denote the event that the random $d$-process $\P(d,n)$ fails to saturate.
 \begin{theorem}\lab{t:unsat} We have
  $$\pr(F)\sim \begin{cases}\displaystyle{\frac{d-1}{ \log n}}\qquad \quad\qquad\mbox{when $dn$ is even}\raisebox{-5mm}{\phantom X}\\  \displaystyle{\frac{(d-1)(d-2)}{ \log^2 n}}\quad\mbox{ when $dn$ is odd}.
  \end{cases}$$
\end{theorem}
Note that in this paper all $\log$s are natural, and we always assume that  $d\geq 2$ is a constant  since all the problems we consider are trivial for $d<2$.
Unless otherwise stated, all asymptotics are as
$n\rightarrow\infty$.

 Theorem \ref{t:unsat} for the case that $dn$ is even confirms a conjecture stated in \cite{RW92}. For calculations and computer simulations of $\pr(F)$ for small values of $d$ and $n$, see the papers of Ruci{\'n}ski~\cite{Ruc90}, and Bali{\'n}ska and Quintas \cite{BQ91}, as well as other papers of Bali{\'n}ska and Quintas cited therein.

Our new approach provides more accurate results on the degree distribution throughout the process.   Let $G_s$ denote the graph in $\P(d,n)$ when it has $s$ edges, and let  $D_i(s)$  denote  the number of   vertices in $G_s$ that have degree $i$, $0\le i\le d$.
The distribution of degrees is given approximately in~\cite[Theorem~3.2]{des}, which we restate as follows.

\begin{theorem}[\cite{des}]\lab{t:old}
For $ 0\le i\le d$,  $\eps>0$, and $s\le dn/2-\eps n$,  we have  $ D_i(s) =nz_i(s)+O(\alpha n)$ with probability $1-O(e^{-n\alpha^3/8})$ for any $\alpha\to 0$, where the functions $z_i$ ($0\le i\le d$) satisfy a certain system of differential equations.
\end{theorem}
The differential equations were specified in that paper and were solved  for $d=2$  in~\cite{RW97} by the authors, and for arbitrary fixed $d$ in~\cite[Section 3.3.3]{des}.  It follows from this solution, by symmetry of the labelling of vertices, that for this range of $s$, the distribution of the degree of any  vertex  of $G_s$ is asymptotically distributed as a truncated Poisson random variable $Z_d(\la)$ (see~\cite[Equation (3.25)]{des}).   To  describe  this distribution, let $Z$  have the Poisson distribution Po$(\la)$ and define  $\pr(Z_d(\la)=i)= \pr(Z=i)$ for $i<d$, and $\pr(Z_d(\la)=d)= \sum_{i\ge d} \pr(Z=i)$. The value of $\la$ is also determined  (see \cite[Formula (3.2)]{des}).

 In the same paper~\cite[Section 5.4]{des}  there is a description of how such results can be extended to cover all $s< dn/2-n^{1-\eps}$, where $\eps$ is a suitable small constant, though no details are given. We extend these results to the following sharper result, applying for a wider range of  $s$, in Section~\ref{s:throughout}.   Let us note at this point, however, that Section~\ref{s:throughout} is not needed  at all for the work in the later sections.

First we need some definitions.
For natural numbers $n$ and real  $x$, let  $\beta_i = \beta_i(n,x) = n\mu^i e^{-\mu}/i!$ where $\mu=x/n$, and define
\bel{ldef}
\ell(n,x)=\frac {dn}2-  \frac12\sum_{i=0}^{d-1}(d-i) \beta_i(n,x).
\ee
We will show in Section~\ref{s:throughout} that for a fixed $n$, $\ell(n,x)$  is strictly increasing in $x$, so it has an inverse function  which we may denote   by $\ell^{-1}_n(s)$.

The following result shows concentration of $D_j$ in the evolving graph of the $d$-process after $s_1$ edges are added, for all $j<d$ provided that the expected number of such vertices tends to infinity.

\begin{theorem}\lab{t:main-degrees}
Fix $d\ge 2$ and $0\le j\le d-1$. Given integer $s =s(n)$ between 1 and $dn/2$, define $x = \ell^{-1}_n(s)$ and  suppose $\beta_j(x)\to\infty$.
 For $w\to\infty$ where $w= o(\beta_j(x)^{1/2})$, we have
$$
\pr\big(|D_{j}( s)-  \beta_j(x)  |\le   w  \beta_j(x)^{1/2}   \big) = 1-O(w^{-2}).
$$
\end{theorem}

For $j=d$  an analogous  result follows by similar arguments. We omit further details, since $D_d(s)=n-\sum_{i<d} D_i(s)$, so a lot of information can be gleaned about $D_d$, especially when $m=\Omega(n)$, immediately from the behaviour of the other degree counts.

The theorem shows that $D_{j}(s)$ is a.a.s.\ inside a window of width growing ever so slightly faster than $ \beta_j(\ell^{-1}_n(s))^{1/2}$ whenever the latter tends to infinity.  This is a narrower window than the previous results using differential equations (as mentioned above). The theorem also immediately extends the range of validity (and precision) of the above conclusion on the truncated Poisson distribution of a vertex degree.

Seierstadt~\cite{S} gave an extension of the main theorem in~\cite{des}, showing that with some mild extra conditions, it follows that the variables concerned converge to a multivariate normal distribution. This is applied in~\cite[Theorem 4]{S}  to the vertex degrees $D_i(s)$ in the random $d$-process, also with the restriction $s\le dn/2-\eps n$. We believe that with some extra work, stronger asymptotic normality results on the joint distribution can be obtained via our approach, since it ties the degree distribution closely to a very simple balls-in-bins model.

  To prove Theorem~\ref{t:unsat}, it is crucial to extend knowledge of the degree distribution even further towards the end of the $d$-process. In particular,   we prove the following three results in Section~\ref{s:survival}.

Define $N=\lfloor dn/2\rfloor$, which is the maximum possible number of edges added in the $d$-process on $n$ vertices.
 For  a real $x$ and an integer $k\ge1$, we set $[x]_k=x(x-1)\cdots(x-k+1)$ for the $k$-fold falling factorial beginning at $x$.
For all $1\le t\le N$ and $0\le j\le d-2$, define
$$
f_j(d,t,n)=\frac{2[d-1]_{d-1-j}t}{\log^{d-1-j}n}.
$$

\begin{theorem}\lab{sharpL7}
For every $j=0,...,d-2$, $k\ge1$, constant $\gamma>0$, and $t=t(n)<\log^\gamma n$ with $t\to\infty$,
$$
 \ex\big([D_{j}( N-t)]_k\big)\sim f_{j}(d,t,n)^k.
$$
\end{theorem}

\begin{theorem}\lab{t:Ddistrib2}
For every $j=0,\ldots, d-2$,    constant $\gamma>0$,   $t=t(n)<\log^\gamma n$ with $t\to\infty$, and $i \ge0 $,
 we have
$$
\pr(D_{j}(N-t)=i)  = e^{-\la} \la^i/i!+o(1)
$$
where $\la=f_{j }(d ,t,n)$.
\end{theorem}

 The distribution of the time that the last vertex of degree $j$ disappears can now be easily determined asymptotically using Theorem~\ref{t:Ddistrib2}.

\begin{corollary}\lab{c:lastvertex} For $j<d-1$, define  $S_j=\max\{s: \delta(G_s)\le j\}$, which is  the last time there is a vertex of degree $j$ in the random $d$-process. Then for each $t>0$ we have
$$
 \pr(S_j< N-t)  = e^{-\la} +o(1)
$$
 where $\la=f_{j }(d ,t,n)$.
  \end{corollary}
It follows that the vertices of degree $j<d-1$ disappear when $t$ is roughly a constant times $\log^{d-1-j}n$.
An approximate version of this result was obtained for the case $d=2$ by Telcs, Wormald and Zhou~\cite{TWZ} using the convergence   property of the  solution  trajectories of the associated differential equations. Since then, the same property has been used in other contexts, and recently, Hofstad~\cite{JH} (inspired by one of us talking on some of the results of the work reported in the present paper) used this very different approach to obtain a slightly stronger version of Corollary~\ref{c:lastvertex}.

 For those interested in other properties of $d$-processes, we mention two particular unsolved problems. In~\cite{RW02}, we showed that the final graph of the $d$-process is a.a.s.\ connected, and  in~\cite[Conjecture~6.2]{regsurvey} the  second author made the specific conjecture that this can be strengthened to being $d$-connected when $dn$ is even. (A more precise conjecture in the case $d=3$ is made in~\cite{RW02}.)
This was also generalised to a much more far-reaching conjecture in~\cite[Conjecture~6.3]{regsurvey}, which is that
 for $dn$  even where   $d\ge 1$ is fixed, the final graph of the $d$-process is
  contiguous to the uniformly random $d$-regular graph. (Two sequences of probability measures on the same sequence of underlying sets   are contiguous if each property a.a.s.\ true for either one of the sequences is also true for the other.) Processes more general than $d$-processes, having unequal upper bounds on vertex degrees, were considered by Molloy, Surya and Warnke~\cite{MSW}, who showed that at least in some cases of these processes, the analogous contiguity property does not hold.

\medskip

\subsection{The bin $d$-process} \lab{s:binprocess}

 In the standard random graph process, the edges are added to $n$ vertices one by one, and the graph existing after $s$ edges are added is  the standard Erd\H os-R\'enyi random graph  $\G(n,s)$. There is an obvious way to embed the $d$-process into this process, by skipping an edge if either of its endvertices already has degree at least $d$. After $s$ edges are attempted to be added, the graph of all edges --- those skipped and not --- is  $\G(n,s)$.  This can be used to give a quite good picture of the $d$-process in the early stages where few edges are skipped: the $d$-process graph  $G_s$  is then very close to the random graph $\G(n, s)$ with the same number of edges. However, by the time the average degree is bounded away from 0, at least a  constant fraction of edges that have arrived have been skipped, and the correspondence between the graphs in the two different processes is so weak as to be useless for most purposes.

 Given the graph $\G(n,s)$ and the order of its   edges arriving, the determination of which edges survive is not simple since whether an edge is skipped depends on how many edges incident with its vertices were skipped before it arrived. So whether an edge is skipped can potentially depend on edges that are quite distant from it, and on the times that they arrived.  Our new approach gives a much tighter coupling than this, of the $d$-process with a simple balls-in-bins model. The vertices in the $d$-process correspond to the bins, and its edges to certain pairs of the balls, with the convenient property that almost all the balls in bins with at most $d$ balls correspond to end vertices of edges in the $d$-process.

 Fix  positive  integers $d$ and $n$. Start with $n$ empty  distinguishable
bins and an infinite supply of  balls. Drop  balls  one by one  and
independently into bins chosen uniformly at random. During the process some balls will be assigned numbers, and some of these  will be deemed
\emph{good}. A bin is {\em unsaturated} if it has less than $d$ good balls in it. Balls
that fall into unsaturated bins are numbered consecutively 1, 2, etc.
Any ball falling into a saturated bin remains unnumbered.

A {\em pair} of balls is defined as the set containing the two balls numbered $i-1$ and $i$, for some even $i\ge 2$. Thus, each ball in the pair is numbered. Call a
pair, as well as each ball in it,  {\em good}  if the balls of this pair fall into distinct bins
and, moreover, no previous pair of balls fell into the same pair of bins.  Otherwise the pair,
and both balls in it, are called {\em bad}. At any time that there is a
numbered ball which is
not yet paired, we call it the {\em waiting} ball. The waiting ball is neither good or bad, even if it
is clear that it will become bad as soon as the pair will be formed. Hence, it does not count in evaluating whether a bin is saturated.

When the bins are nearly all saturated, most balls will be unnumbered because they  keep dropping into saturated bins.  The process terminates when no more good pairs can be added. Note that at each point in time, each ball is classified as exactly one of good, bad, waiting, or unnumbered.

This random process  is called the {\em bin $d$-process} on $n$ bins. We associate a graph process
with this bin $d$-process as follows.
 The vertices are the bins, while the good pairs determine edges
joining the bins into which these balls have fallen.   Let $G'_s$ be the graph obtained after $s$ edges are added. Then $\B(d,n):=(G'_0,G'_1,\dots)$  is the {\em
bin-graph $d$-process} on $n$ bins. Note that the degree of each vertex in $G'_s$
is determined by the number of good balls in the corresponding bin.

The following lemma, though quite obvious, lies at the heart of our  method.

\begin{lemma}\lab{l:equiv} For all positive integers $d$ and $n$,
the bin-graph $d$-process $\B(d,n)$ on $n$ bins has exactly the same distribution as the  graph $d$-process
$ \P(d,n)$.
\end{lemma}

\proof Given the first $ s-1$ pairs of numbered balls, each of the next two numbered balls is equally
likely to fall into any  unsaturated bin. Hence, conditional on being good, these two balls are equally likely to lie in
any pair of distinct  unsaturated  bins. So the $s$th pair of good balls has exactly the same distribution,
given that   $G'_{s-1}$ is equal to a specific graph $H$, as the $ s$th edge of
$\G(d,n)$,  given $G_{s-1}=H$. Hence the distributions of the two random  graphs  remain
identical throughout the two processes. \qed

Based on the above lemma, we often in this paper view the  $d$-process $\G(d,n)$ as the bin-graph $d$-process derived from the bin $d$-process.  (In particular, from now on we use the notation $G_s$ for $G_s'$.)  Hence we may speak of the balls and bins  when studying the former.
 This equivalence enables us to use straightforward properties of the standard balls-in-bins random
allocation model to derive a lot of information on the degree distribution, which we will use to estimate the asymptotic probability of  non-saturation of the graph $d$-process $\P(d,n)$ as stated in Theorem \ref{t:unsat}.

 In our proofs it will be convenient to analyse the process in terms of the residual number of (potential)  good balls which still may occur before the process ends. Let $\widetilde X(m)$ be the number of good balls among the first $m$ balls dropped in (note that $\widetilde X(m)$ is always even).
At each point $m$ in time (here time is measured by the number of balls dropped), we define the \emph{deficit}
 $T(m)$ of the bin $d$-process to be  the maximum possible number of good pairs of balls minus the  number of good pairs of balls present at time $m$, that is,
\bel{Tdef}
T(m)=  N - \widetilde X(m)/2.
\ee

\section{Vertex degree  versus number of balls in a bin}\lab{s:distrib}

In this section we will point out some easy connections between the number of
vertices of a given degree $i$ and the number of bins containing exactly $i$  balls. Tighter estimates
will be made later, as they are needed for close examination of the later part of the bin process. Our first step is to prove a very useful lemma about the vertex degrees in the standard (graph) $d$-process.

 Note that when the deficit in the bin $d$-process reaches $t$, the number of good pairs of balls is exactly $N-t$, and hence this is also the number of edges  $s$  that have been added in the corresponding bin-graph $d$-process. We call $N-s$ the {\em deficit} of the $d$-process, since it is the deficit in every corresponding bin $d$-process. The deficit of the graph process is simply the maximum number of edges that can possibly be added later in the process, based on the degree bounds.

 It follows that for a fixed deficit $t$, the number of good balls  $\widetilde X(m)$  is independent of $m$ -- the actual number of balls dropped in so far, so we may denote it by $X(t)$. Similarly,  regardless of $m$, let $U(t)$ be the number of unsaturated bins present at deficit $t$.
Then, by definition of  deficit,
\bel{Ubounds}
2t+1\ge U(t)\ge 2t/d.
\ee
Indeed, denoting by $g_j$, $j=1,\dots U(t)$, the numbers of good balls in the unsaturated bins at deficit $t$, we have
$$2t\le dn-X(t)=\sum_j(d-g_j)\le dU(t)$$
and
$$2t\ge dn-1-X(t)=\sum_j(d-g_j)-1\ge U(t)-1.$$

 \subsection{Survival of unsaturated vertices}
 Our next result will allow us  later to study the process at very small deficits. It will be also instrumental in estimating the number of bad balls.
\begin{lemma}[Survival Lemma]\lab{l:survival}
Let $k$  and $a$ be two fixed integers and $v_1,\ldots , v_k$ vertices, and condition on a state with deficit $t_0$ in the $d$-process in which $v_1,\ldots , v_k$ are all unsaturated.   Then, provided that $ t_0>t\to\infty$, the probability  that, by the time the deficit reaches $t$, the total of the degrees of these vertices has increased by $a$ and they all have remained unsaturated, is
$$
O(1)\big(\log (t_0/t)\big)^{a} \left(\frac{t}{t_0}\right)^k.
$$
In particular, the probability they all remain unsaturated until the deficit reaches $t$ is
$$
O(1)\big(\log (t_0/t)\big)^{dk} \left(\frac{t}{t_0}\right)^k.
$$
\end{lemma}
\proof
 The event $A$ that the   vertices $v_1,\ldots , v_k$ gain total degree $a$ by  deficit   $t$  and remain unsaturated breaks into subevents $A(M)$ where $M$ is a matrix $(r_{i,j})_{1\le i\le k, 1\le j\le d-1}$ and $r_{i,1},\ldots, r_{i,d-1}$ denote the  deficits from $t_0$ down to $t+1$ at  which the next edge added is incident with  $v_i$.  (We may set  some   of these variables equal to infinity to denote that such increases do not occur.) Note that a vertex degree can only increase by 1 in a given step  since loops are forbidden in the $d$-process.

Also, for $t'\ge  t$,  denote by  $A_{t'}(M)$ the event that the changes in the degrees of the $v_i$ are consistent with the prescription of the event $A(M)$ for all the steps up to and including the one in which   the deficit reaches   $t'$.
For any $t'>t$ that is not equal to any of the $r_{i,j}$,  conditional on the history of the process up to  deficit   $t'$, provided the event  $A_{t'}(M)$ holds, the probability of
$A_{t'-1}(M)$ is equal to
$$
1-\frac{2k}{U(t' )+O(1)}
$$
since there are $kU(t')+O(1)$ edges that would increase the degree of one of the $v_i$, and ${U(t')\choose 2}-O(U(t'))$ edges that can be chosen in this step.
Similarly, if $t'$ is equal to any of the $r_{i,j}$ (but not more than one), we get an expression
$$
 \frac{2}{U(t')+O(1)}.
$$
If two of the $r_{i,j}$ are equal to $t'$, the expression is at most
$$
\frac{O(1)}{(U(t')+O(1))^2},
$$
and of course it is not possible that the degrees of three vertices are changed in one step.
Now using the upper and lower bounds in~\eqn{Ubounds},   the  expressions displayed above are bounded above by
$$
1-\frac{2k}{2t'+O(1)},\quad
 \frac{O(1)}{t'}\quad {\rm and}\quad
 \frac{O(1)}{(t')^2}
$$
respectively.
These bounds are independent of $U(t')$. Hence  there is no need to condition on the history any more, and we deduce using the chain formula for conditional probabilities that
$$
\pr(A_t(M)=
O(1)\bigg(\prod_{i,j\,:\,r_{i,j}<\infty}\frac{1}{r_{i,j}} \bigg) \prod_{t'=t+1 }^{t_0}\left(1-\frac{2k}{2t'+O(1)}\right),
$$
where we note that the extra factors introduced into the second  product, for those values of $t'$ equal to one of the $r_{i,j}$, do not affect the expression.

Given $M$, let $c_i$ denote the total change  in degree  of the vertex $v_i$ ($1\le i\le k$) as the deficit changes from  $t_0$ to $t$, i.e.\ the largest $j$ for which $r_{i,j}<\infty$.
Summing over all possibilities for the numbers $r_{i,j}$ in the matrix $M$ that are feasible given the $c_i$,  using $1-x\le e^{-x}$  and recalling that $t\to\infty$  (so the denominator $2t'+O(1)$ is not zero) we obtain
\bean
\pr(A)&=& O(1) \exp\bigg(-\sum_{t'=t+1 }^{t_0}  \frac{ k}{ t'}\ \bigg)\sum_{M}
\bigg(\prod_{i=1}^{k}\prod_{j=1}^{c_i} \frac{1}{r_{i,j}}\bigg)\\
&=& O(1) \exp\bigg(-k\sum_{t'=t+1 }^{t_0}  \frac{1}{ t'}\ \bigg)
\sum\sum\cdots \sum
\prod_{i=1}^{k}\prod_{j=1}^{c_i} \frac{1}{r_{i,j}}
\eean
where the $i$-th sum is over  $t <r_{i,1}<\cdots<r_{i,c_i}\le t_0$.  We can easily bound this above by relaxing the ordering on the $r_{i,j}$ and permitting them to be equal. The result is an expanded version of a power of a single sum, giving
\bean
\pr(A)&=& O(1) \exp\bigg(-k\sum_{t'=t+1 }^{t_0}  \frac{1}{ t'}\bigg)
 \biggl(\sum_{r=t+1 }^{t_0}\frac{1}{r }  \biggr)^{c_1 +\cdots +c_k}.
 \eean
Each summation is $O(1) + \log (t_0/t)$ and is also  $O(\log (t_0/t))$. As the total change in degrees of the $v_i$ is  $a = c_1 +\cdots +c_k$,   the result follows.
 \qed

 \subsection{Bad balls}

 Consider the bin $d$-process on $n$ bins.  Let ${\cal S}_m$  denote the state of the    process after $m$ balls
have been dropped in, where one might have $m=m(n)$.    For $i=1,\dots,d-1$, define the following statistics of ${\cal S}_m$.  (Here, and occasionally later, we suppress the dependence on $m$.)

\begin{itemize}
\item $Y_i $: the number of bins containing $i$ balls (good, bad, waiting, or unnumbered);
\item $\DD$: the number of bins containing $i$ good balls, i.e.\ the number of vertices of
degree $i$ in the corresponding graph;
\item $B$:  the  number of bad  balls.
\end{itemize}

\begin{claim}\lab{determ} At any time  in the process,
\bel{Dapprox} \mbox{ $|Y_i-\ \DD |\le  B+1$ for $0\le i\le d-1$,}\ee
\end{claim}

\proof
For every bin containing no waiting or bad ball, the degree of the corresponding vertex equals the number of balls in the bin. Hence, estimating  $\DD$  by $Y_i$ gives an error no larger than the number of bins containing a bad or waiting ball, of which there can be at most $B+1$. \qed
 \smallskip

By definition of the process, the distribution of the numbers of balls in the bins --- both numbered and unnumbered in total ---  is exactly
multinomial with parameters $m$ and $(1/n,\ldots,1/n)$.
In the next section, we will show that the random variable $Y_i$ is   concentrated, provided its expected value is large, and our main problem will be to establish a link between $Y_i$ and $\DD$.
In view of~\eqn{Dapprox}  this can be done by bounding the number $B=B(m)$ of bad balls among the first $m$ balls dropped.
We also consider the  \emph{total} number of bad balls occurring by  the end of the process, which we denote by $\Bf$. In the lemma below, we  estimate   its  $k$-th factorial moment, which we will find useful.

 \begin{lemma}\lab{badballs} For every $k\ge1$,  $\ex (\Bfk) = O(\log^k n)$.
\end{lemma}
\proof
 For $t\ge0$, let $W_t$ be the number of bad pairs of balls   created while the deficit is $t$.
Then
  $$\Bf=2\sum_{t}W_t\quad\mbox{ and so}\quad\ex \Bfk=2^k\ex\left(\sum_t W_t\right)^k.$$
  Our  goal  is to bound from above the last expression. By the multinomial theorem and linearity of expectation,
\bel{summation}
\ex \left(\sum_t W_t\right)^k\le C \sum_{s=1}^k\sum_{k_1+\cdots + k_s=k}\ \sum_{0\le t_1<\cdots<t_s\le t_0}\ex\left(\prod_{j=1}^{s}W_{t_j}^{k_j}\right),
\ee
where the $k_i$ define an ordered partition of $k$, and $C$ is an upper bound on the multinomial coefficients.  (Here $C$  is constant, since $k$ is fixed.)

We will show by induction on $s$ that for all $0<t_1<\cdots<t_s$
\begin{equation}\lab{ind}
\ex(W_{t_1}^{k_1}\cdots W_{t_s}^{k_s})=O((t_1\cdots t_s)^{-1}).
\end{equation}
 For  $t=t_1>d^2$  note that $U=U(t)>d$ by~\eqn{Ubounds}.
   When a waiting ball is sitting in a bin,  this ball and the next  numbered ball will turn bad exactly when the latter has landed in the same bin as the waiting one or in one of the adjacent bins. This has probability at most $d/U$, since (at any time) each unsaturated bin is adjacent, in the corresponding graph, to at most
 $d-1$ other bins.
  Thus, $W_t$ (conditioned on the state of the process ${\cal S}(m)$ for  $m$ when the deficit first reaches $t$)
  is stochastically dominated by a geometrically distributed random variable $W$   where $\pr(W=i) = p(1-p)^i$ with $p=1-q=1-d/U$ and thus  with mean
 $$q/p=d/(U-d)<\frac{d^2}{2t-d^2}<\frac{d^2}t=O(1/t)$$
  by~\eqn{Ubounds}.  Moreover, for any fixed $r$,
$$
\ex(W_{t}^r)\le\ex(W^r)= \sum_{i\ge 1} pq^i i^r \le \sum_{i\ge 1} pq^i r!{i+r-1\choose r}=  r!pq(1-q)^{-(r+1)} =r!p^{-r}q$$
which is $O(q)=O(1/U) =O(1/t)$ for fixed $r$.

  In the complementary case, when $t\le d^2$, we instead use the fact that there is  always at least one potential good pair remaining, so we can use the argument above  with $W$ replaced by a geometrically distributed random variable  $W'$ with
$q'=1- p' =1- \binom U2^{-1}$, where $U=O(1)$
by~\eqn{Ubounds}. So, $(p')^{-1}=O(1)$ and, consequently,
$$
\ex(W_{t}^r)\le\ex (W')^r\le r!(p')^{-r}q' =O(1)=O(1/t),
$$
as before.

Thus, in each case, we have $\ex W_{t_1}^{k_1}  = O(1/t_1)$
which establishes \eqref{ind}   for $s=1$. Now, assume it is true for $s-1$, $s\ge2$.  Then, by the induction hypothesis, setting $\hat W=W_{t_1}^{k_1}\cdots W_{t_{s-1}}^{k_{s-1}}$,
 $$\ex(W_{t_1}^{k_1}\cdots W_{t_s}^{k_s})=\ex(\hat W)\ex(W_{t_s}^{k_s}|\hat W)=O\left(\frac1{t_1\cdots t_{s-1}}\right)\ex(W_{t_s}^{k_s}|\hat W).$$
Since the above argument about $W_{t_1}$ also applies to  $W_{t_s}$ conditioned on the state of the process ${\cal S}(m)$ for  $m$ when the deficit first reaches $t_s$, we similarly obtain $\ex(W_{t_s}^{k_s}|\hat W)=O(1/t_s)$, so
 $$\ex(W_{t_1}^{k_1}\cdots W_{t_s}^{k_s})=O\left(\frac1{t_1\cdots t_{s}}\right),$$
 as desired.

 In view of \eqn{summation}, we finish off by summing over  all $s$ and all $t_1,\dots, t_k$. Thus,
  $$\ex \left(\sum_t W_t\right)^k=\sum_{s=1}^k\sum_{t_1}^\infty\cdots\sum_{t_s}^\infty\frac{O(1)}{t_1\cdots t_s}=\sum_{s=1}^kO(\log^sn)=O(\log^kn).$$
  The proof is complete. \qed

By combining the Survival Lemma  (Lemma \ref{l:survival}) with ideas from the proof of Lemma \ref{badballs}, we can obtain a more accurate estimate of number of bad balls that are relevant to our investigations.
Let $\widetilde B(t)$ denote the number of bad balls in {\em unsaturated}  bins  at  the last step in the bin process when the deficit is $t$.
Note that the proof of Claim \ref{determ} yields, in fact, that
\bel{Cl22t}
|Y_i-\ \DD |\le  \widetilde B(t)+1
\ee
for $0\le i\le d-1$.

\begin{lemma}\lab{c:betterbad}
Uniformly for all $t \ge 0$ we have $\ex \widetilde B(t)=O(1)$.
\end{lemma}
\proof
From~\eqn{ind} with  $s=1$, $k_1=1$,  and $t_1=t_0$  we have that for any $t_0\ge 0$, the expected number of vertices containing  bad balls created when the deficit is $t_0$ is $O(1/t_0)$.   Conditioning on the graph at the last time in the ball process that the deficit is $t_0$,  the Survival Lemma gives that the probability that any such vertex survives to be unsaturated at deficit $t$ is $O(1)\big(\log (t_0/t)\big)^{d-1} \left(\frac{t}{t_0 }\right) $. Hence, using  linearity of expectation, we see that the expected number of   vertices unsaturated at deficit $t$ and containing a bad ball added when the deficit was $t_0$ is
$$
O(1)\big(\log (t_0/t)\big)^{d-1} \left(\frac{t}{t_0^2}\right) .
$$
Now $\ex \widetilde B(t)$ is the sum of this quantity over all $t_0\ge t$. Hence,  substituting $u=t_0/t$,   we have
\begin{align*}
\ex \widetilde B(t) &= O(1) \int_{t}^\infty  \big(\log (t_0/t)\big)^{d-1} \left(\frac{t}{t_0^2}\right)dt_0 \\
& = O(1) \int_{1}^\infty \frac{\log^{d-1} u}{u^2} du  = O(1),
\end{align*}
since the integral evaluates to exactly $(d-1)!$. \qed

 It follows by Markov's inequality that uniformly for all $b=b(n)>0$ and $0\le t\le N$         we have
\bel{secondk}
\pr(\widetilde B(t)\ge b ) = O(b^{-1}).
\ee

With some work, one can combine the ideas in the proofs of Lemmas~\ref{badballs} and~\ref{c:betterbad} to obtain strong bounds on the higher moments of $\widetilde B(t)$. This would result in better bounds on the degree distribution of vertices in the $d$-process, but the extra effort is not warranted for our main objectives in this paper.


\section{ Concentration of  $Y_i$  in the bin $d$-process}
Consider the bin $d$-process on $n$ bins   at state ${\cal S}(m)$, that is, a uniformly random allocation of $m$ balls into $n$ bins.
Let $A_j$ be the number of balls in bin $j$, $j=1,\dots,n$. Recall that $Y_i=| \{j:A_j=i\}|$.  By linearity of expectation  (considering  the probability that the first bin receives $i$ balls), provided $m=o(n^2)$ we have
\bel{EY}
\ex Y_i =n\frac{\binom mi(n-1)^{m-i}}{n^m}\sim n\mu^i e^{-\mu}/i!=  \beta_i(n,m)
\ee
as defined in Section~\ref{s:intro}, with $\mu=m/n$.
Similarly, by considering the event that two bins receive $i$ balls each,
$$
\ex [Y_i]_2 =n(n-1)\frac{\binom {m}{i}\binom {m-i}{i}(n-2)^{m-2i}}{n^m}
\sim n^2\mu^{2i} e^{-2\mu}/i!^2= \beta_i^2.
$$
In fact, with  slightly more careful analysis, we will show shortly that
\bel{EY2}
 \ex  Y_i  = \beta_i+O(1)  \quad \mbox{ and } \quad \ex [Y_i]_2 = \beta_i^2 + O(\beta_i+1).
\ee
 Indeed,
$$\ex Y_i=\beta_i\frac{[m]_i}{m^i}e^{m/n}\left(1-\frac1n\right)^{m-i}=\beta_i\left(1+O\left(\frac1m+\frac1n+\tfrac m{n^2}\right)\right)$$
and similarly
$$\ex[Y_i]_2=(\ex Y_i)^2\left(1+O\left(\frac1m+\frac1n+\frac m{n^2}\right)\right).$$
Now, if $\mu=m/n=O(1)$, then $m/n^2=O(1/n)=O(1/m)=O(1/\ex Y_i)$, as $\ex Y_i\le m$ trivially. On the other hand, if $\mu\to\infty$, then $1/m=o(1/n)=o(m/n^2)=o(1/\beta_i)$, because $\mu^{i+1}=o(e^{\mu})$. This yields \eqref{EY2}, since  $\ex Y_i=\beta_i(1+O(1/\ex Y_i))=\beta_i+O(1)$ and so
$$
\ex[Y_i]_2=(\ex Y_i)^2+O(\ex Y_i)=\left(\beta_i+O(1)\right)^2+O(\ex Y_i)=\beta_i^2+O(\beta_i +1).
$$
It follows from \eqref{EY2} that
\bel{VarY}
{\bf Var\,} Y_i = \ex [Y_i]_2+\ex Y_i-(\ex Y_i)^2   = O(\beta_i+1).
\ee
provided $m=o(n^2)$.

This  gives us our first approximation of $Y_i$ by $\beta_i$.
\begin{lemma}\lab{l:easy}
 For any $ \alpha=\alpha(n)  >0$, if  $i$ is fixed,  $m=o(n^2)$, and $\beta_i=\Omega(1)$,   we have
 $$\pr( |Y_i -\beta_i|\ge  \alpha\sqrt{\beta_i})=O(1/\alpha^2).$$
\end{lemma}
\proof We may assume that $\alpha\to\infty$ as otherwise the claimed bound is trivial.
By~\eqref{EY2}, there is a constant $C>0$ such that $|\ex Y_i-\beta_i|<C$. By \eqref{VarY},    ${\bf Var Y_i\,}= O(\beta_i+1)$, and so, by
Chebyshev's inequality, for sufficiently large $n$, this gives
\begin{align*}& \pr(|Y_i -\beta_i|\ge \alpha\sqrt{\beta_i})\le\pr(|Y_i -\ex Y_i|+C\ge \alpha\sqrt{\beta_i})\le \pr\big(|Y_i -\ex Y_i|\ge \frac{\alpha}2\sqrt{\beta_i} \big)\\
 &\le \frac{{\bf Var\,} Y_i}{ \alpha^2\beta_i/4}
  =O(1/\alpha^2). \qed
 \end{align*}

We also require some sharper concentration of $Y_i$, for use at certain times in the process in an argument that transfers the knowledge to later times in the process when there is less concentration.  For this we  apply  a standard approach involving independent Poisson random variables as follows.

 Clearly, for any sequence of
nonnegative integers $a_1,\dots,a_n$ with $\sum_{j=1}^na_j=m$,
$$\pr(A_j=a_j, j=1,\dots,n)=\frac{m!}{\prod_{j=1}^na_j!}n^{-m}.$$
From this it is easy to verify the well known fact that
$$\pr(A_j=a_j, j=1,\dots,n)=\pr\bigg(Z_j=a_j, j=1,\dots,n \bigg|\sum_{j=1}^nZ_j=m\bigg),$$
where $Z_1,\dots,Z_n$ are iid random variables with Poisson distribution ${\rm Po}(\mu)$ (and this would still hold with ${\rm Po}(\lambda)$  for any $\lambda>0$, but $\mu$ suits our purposes).

 We  observe that $\sum_{j=1}^nZ_j$ has distribution ${\rm Po}(m)$ and thus
$$
\pr\bigg(\sum_{j=1}^nZ_j=m\bigg)=\frac{e^{-m}m^n}{m!}=\Theta(m^{-1/2}).
$$
This equation has the following consequence.
 Let $\pr_A$ and $\pr_Z$ be, respectively, probability measures defined by the sequences of random variables $(A_1,\dots,A_n)$ and $(Z_1,\dots,Z_n)$. Then, for every event $E$,
\begin{equation}\lab{conclude}
\pr_A(E)=O(\sqrt m)\pr_Z(E).
\end{equation}

We can thus deduce distributional results on $Y_i$, under the probability function $\pr_A$, by first obtaining results on $\tilde Y_i $, where $\tilde Y_i =| \{j:Z_j=i\}|$, under the probability function $\pr_Z$.  Observe that $\tilde Y_i$ has a binomial distribution with $\ex\tilde Y_i=\beta_i$.
The next result gives two special cases of the concentration of $Y_i$ near its approximate expected value $\beta_i$ that are useful for us. The proof is a simple application of~\cite[ Corollary 2.3, (2.9)]{JLR}: for any  $\alpha\le  \tfrac32 \beta_i$,
\bel{sharpconc}
\pr(|\tilde Y_i-\beta_i|\ge \alpha) 2 \exp (-\alpha^2/ 3\beta_i).
\ee
We only consider $i\le d-1$ since  the remaining  bins, with at least $d$ balls,  are not examined in any detail in our arguments.
\begin{lemma}\lab{YtoYtilde} Let $i\ge 0$ be fixed. We have
 \begin{enumerate}
\item[(i)] if $\beta_i/\log n\to\infty$, then
$$\pr(| Y_i- \beta_i|\ge 3\sqrt{\beta_i\log n})=o(n^{-2}).$$

\item[(ii)] if, for some $\eps>0$, $\beta_i>\log^{2+\eps}n$, then
 for every $K>0$,
 $$
 \pr(| Y_i- \beta_i|\ge \beta_i^{3/4})= o(\beta_i^{-K}).
 $$
\end{enumerate}
\end{lemma}
\proof Both statements will follow from a more general one which is an immediate corollary of \eqref{conclude} and \eqref{sharpconc}.
 To transfer \eqref{sharpconc} to tail bounds on $Y_i$, note that the assumption $\beta_i/\log n\to\infty$ implies that $m=O(n\log n)$. Indeed, $\beta_i=ne^{-m/n}(m/n)^i/i!$ which, by setting $m=an\log n$, $a=a(n)$, becomes $\beta_i=a^in^{1-a}\log^i n/i!$. For this quantity to tend to $\infty$,  it is necessary that  $a=O(1)$.

 With this bound on $m$, equations \eqref{conclude} and~\eqn{sharpconc} above now yield
  \bel{gen}\pr(|Y_i-\beta_i|\ge\alpha)=O(\sqrt{n \log n})\pr(|\tilde Y_i-\beta_i|\ge\alpha) =O(\sqrt{n \log n})  \exp (-\alpha^2/ 3\beta_i),
  \ee
  which implies both concentration results   (i) and (ii). Indeed, with $\alpha=3\sqrt{\beta_i\log n}$, \eqref{gen} yields in (i) the bound $O\left(n^{-2.5}\sqrt{\log n}\right)=o(n^{-2})$, while setting $\alpha=\beta_i^{3/4}$, the probability in (ii) can be bounded by
  \begin{align*}
  &O(\sqrt{n \log n})\exp(-\tfrac13\sqrt{\beta_i})=O(\sqrt{n \log n})\exp(-\tfrac13\log^{1+\eps/2}n)\\&=O(\sqrt{n \log n})\left(\frac1n\right)^{\tfrac13\log^{\eps/2}n}=o(n^{-K})=o(\beta_i^{-K}),
  \end{align*}
  because $\beta_i\le n$ by definition.
  \qed

The main results in this paper require estimates of $\ex Y_i$, via $\beta_i$, at times very near the end of the process, when the number of unsaturated bins is at most a power of $\log n$.
Recall that $\beta_i= ne^{-\mu}\mu^i/i!$ and  $[x]_k$  denotes the falling factorial $x(x-1)\cdots(x-k+1)$. The following relations are useful to isolate.

\begin{claim}\lab{c:beta}
 If   $0\le i \le d-1 $ then
\bel{tildeYasy}
 \beta_{i}= \beta_{d-1} [d-1]_{d-1-i}/\mu^{d-1-i}.
\ee
Moreover, if  $\mu=\Omega(1)$ and $\beta_{d-1}=\log^{O(1)} n $, then $\mu\sim \log n$ and
\bel{CYasy}
\beta_i\sim  \beta_{d-1} [d-1]_{d-1-i}/\log^{d-1-i}n. \ee
\end{claim}

\proof
The first assertion is immediate from the definition of $\beta_i$ in~\eqn{EY}. For the second,  observe that if   $\mu=\Omega(1)$ and   $\beta_{d-1}=\log^{O(1)} n $ then  $\mu\sim \log n$. Consequently, ~\eqn{tildeYasy} gives $ \beta_{i} \sim  \beta_{d-1} [d-1]_{d-1-i}/\log^{d-1-i}n$,  which is \eqref{CYasy}. \qed

\section{Degree distribution throughout the process} \lab{s:throughout}

This section considers degree counts $D_j$ in the $d$-process for $s=s(n)$ such that $D_j(s)\to\infty$  in expectation.  The argument here has issues, of applying the bin model, that are in common with the work in the later sections, which deals with   the very late steps in the process  where $D_j$ is sometimes likely to be bounded.
In this section, we suppress the dependence of all variables on $n$, since it only plays a minor role.
In particular, we will write $\beta_i(m)$ instead of $\beta_i(n,m)$.

We restate Theorem \ref{t:main-degrees} as follows with $s$ and $x$ replaced by $s_1$ and $m_1$ for convenience of the proof.

 \newtheorem*{thm:1point4}{Theorem \ref{t:main-degrees}}
\begin{thm:1point4}
 Fix $d\ge 2$ and $0\le j\le d-1$. Given any $s_1=s_1(n)$ between 1 and $dn/2$, suppose $\beta_j(m_1)\to\infty$ where $m_1 = \ell^{-1}_n(s_1)$.
 For $w\to\infty$ where $w= o(\beta_j(m_1)^{1/2})$, we have
\bel{Dconc}
\pr\big(|D_{j}( s_1)-  \beta_j(m_1)  |\le   w  \beta_j(m_1)^{1/2}   \big) = 1-O(w^{-2}).
\ee
\end{thm:1point4}

In the rest of this section we prove this theorem. (We also  show that the   inverse function   of $\ell$ exists.)
Define
\begin{align}\lab{Ldef}
L(m) &=\frac12 \sum_{ i=0}^{m} \min\{i,d\}Y_i(m))=
\frac12 \sum_{i=0}^{d-1} iY_i(m))+ \frac d2\bigg(n-\sum_{i=0}^{d-1} Y_i(m)\bigg)
\nonumber\\
 &=  \frac{dn}{2}-  \frac12\sum_{i=0}^{d-1}(d-i) Y_i(m).
\end{align}
 Thus, $L(m)$ is one half of the current number  of balls if we ignore balls that had arrived  in a bin which already contained at least $d$ balls. Note  that one can write
$$L(m)=\frac12\sum_{ i=1}^d\sum_{j=i}^mY_j(m)$$
which shows that $L(m)$ is a non-decreasing function of $m$, because   $Y_{\ge i}(m):=\sum_{j=i}^mY_j(m)$ is   non-decreasing for every $i$.

 Let  $\widetilde S(m)$   be the (random) number of edges added to the graph in the bin $d$-process
 due to the first $m$ balls dropped.
  It follows that $\widetilde S(m)\le L(m)$, as edges correspond to the pairs of good  balls placed in the bins, and there are at most $d$ good balls in any bin.
  Further, let $\widetilde M(s)$   be the number of balls  in the bins  at the last step where the number of edges is $s$ (and hence the deficit is $t=N-s$).  By this definition, the waiting ball is necessarily present  among the $\widetilde M(s)$ balls. Recall from Section \ref{s:distrib} that $\widetilde B(t)$ denotes the number of bad balls in unsaturated bins at the last step in the bin process
when the deficit is $t$. We have
\bel{L-to-s}
s=L(\widetilde M(s))-\widetilde B(N-s)/2 -\hat \delta/2
\ee
where $\hat \delta$ is 1 if the waiting ball is present (which occurs in all but the last step) and 0 otherwise.
 Since  $\widetilde B$ is expected to be small according to Lemma~\ref{c:betterbad},   $L(m)$ should be close to the number  $\widetilde S(m)$  of pairs of good balls  (though the fact that the lemma is in terms of one particular value of $t=N-s$, to which many different values of $m$ can be relevant, complicates the discussion shortly).

  Recall the definition  of $\beta_i$ in  \eqref{EY}.  As $\beta_i$ is a good approximation for   $\ex Y_i$ in view of~\eqn{EY2}, we will estimate  $L(m)$, and through it $\widetilde S(m)$,
 by
$$
\ell(m)=\frac {dn}2-  \frac12\sum_{i=0}^{d-1}(d-i) \beta_i(m)
$$
as defined in~\eqn{ldef}.
 Here and in the argument which follows we will consider $m$ to be non-negative real, and we define $\tau=\tau(m) = \sum_{i=0}^{d-1} \beta_i (m)  $.

 We use the following guarantee that $L$ is well approximated by $\ell$ as a function of $m$.
\begin{lemma}\lab{l:l}
For $  \alpha>1  $ and any $m >0$ with $m=o(n^2)$,
$$
\pr\big(|L(m)-\ell(m)|>\alpha \tau(m)^{1/2}\big)=O(1/\alpha^2).
$$
\end{lemma}
\proof
The crucial point is that
$$
|L(m)-\ell(m)|\le \frac{d}{2}\sum_{i=0}^{d-1} | Y_i(m) -\beta_i(m)|.
$$
 So $|L(m)-\ell(m)|>\alpha \tau(m)^{1/2}$ only if for some $i$
\bel{Y-beta}
|Y_i(m)-\beta_i({m})|> 2\alpha \tau(m)^{1/2}/d^2 \ge 2\alpha \beta_i(m)^{1/2}/d^2.
\ee

The probability of this event for a given $i<d$ is $O(1/\alpha^2)$ by   Lemma~\ref{l:easy}  (with $\alpha$ replaced by $2\alpha/d^2$).
The union bound (over $0\le i\le d-1$) now implies that with probability $1-O(\alpha^{-2})$, none of the $d$ events~\eqn{Y-beta} holds, which gives the lemma.
\qed

 The lemma shows that the typical variation in $L(m)$ is   $O(\tau(m)^{1/2})$.
Note that $\tau  = \Theta(\max \{ \beta_0  , \beta_{d-1}  \})$, and that, considering~\eqn{tildeYasy}, $\beta_0$ attains this maximum   near the start of the process, and $\beta_{d-1}$ near the end.

The main complication in  the proof of Theorem \ref{t:main-degrees}  is due to the fact that~\eqn{L-to-s}, which links $s$ and $L$, is based on a given value of $s$, the number of edges in the graph, whereas Lemma~\ref{l:l}, linking $L$ and $\ell$, is based on  a given value of $m$. We will deal with this problem by defining two values of $m$, called $m_0$ and $m_2$ below, such that a given number $s_1$ of edges is highly likely to be attained in the graph when the number of balls is between $m_0$ and $m_2$.

 To proceed we need to know that function $\ell(m)$ is invertible, and also to bound some derivatives. To this end, note  that
  \bel{betaderiv}
   \beta'_i(m)  = (i/\mu-1)\beta_i(m)/n=  (\beta_{i-1}(m)-\beta_i(m))/n
\ee
 (with $\beta_{-1}$ interpreted as 0), and so, for  real  $m>0$,
\bel{lderiv}
\ell'(m) =
 \frac {1}{2n}\sum_{i=0}^{d-1}(d-i) (\beta_{i}(m)-\beta_{i-1}(m)) = \frac{\tau(m)}{2n}>0.
\ee
 Thus $\ell(m)$  is strictly increasing in $m$, as required to show that the inverse of $\ell$ exists.
 Of course  $\ell^{-1}_n(s)$ is usually not an integer, but we may use it as  an estimate of $\widetilde S(m)$.

 To avoid potential oscillation of the random variable $Y_i(m)$  in value as $m$ increases, we instead  study $Y_{\ge i}(m):=\sum_{j=i}^mY_j(m)$ and $Y_{\le i}(m):=\sum_{j=0}^iY_j(m)$, which are (weakly) monotonically increasing and decreasing respectively.  After this, we estimate $Y_i$ via $Y_{\le i}(m)-  Y_{\le i-1}(m)$ or $Y_{\ge i}(m)-  Y_{\ge i+1}(m)$.  We need to   consider  two cases.    The first occurs when $Y_i$ is essentially increasing  in $i$, which happens when $m$ is large, and the second when $m$ is small.   If $m=\Theta(n)$, the values for different $i$ will be of comparable size and the argument in either case will apply.

Recall that $m_1=\ell^{-1}(s_1)$, as defined  in the  (re)statement of Theorem~\ref{t:main-degrees}.

\medskip

\noindent
{\bf Case 1: $m_1\ge n$.\/}
\smallskip

Choose a  function $w=w(n)\to \infty$   with $w=o (\tau(m_1)^{1/2})$. Define
$$m_0= m_1 -wn/\tau(m_1)^{1/2}\quad\mbox{and}\quad m_2= m_1 +wn/\tau(m_1)^{1/2}.$$

 We will next show that it is highly likely that $  m_0 <\widetilde M(s_1)<m_2$, (i.e.\ the number of balls dropped when the graph has $s_1$ edges   is between $m_0$ and $m_2$), and that the random variables $\Db_i(m)$ are sufficiently concentrated for our purposes whenever $m_0\le m\le m_2$.

For $m_0\le m\le m_2$,   as $\mu=\Omega(1)$, we have from the first equation in \eqn{betaderiv} that $ \tau'(m) = O(\tau(m)/n)$. Solving the differential equation $ \tau'(m)=c\tau(m)/n $  we obtain,  for $m_0\le m\le m_2$ and some $c,C>0$, $\tau(m)=e^{cm/n+C}$.  It follows that $\tau(m )=\tau(m_1)e^{c (m -m_1)/n }$.
Thus, noting that $|m_1-m|\le wn/\tau(m_1)^{1/2}$, we have
\bel{pre-conc}
\tau(m) = \tau(m_1)e^{O\left( |m_1-m| /n\right)}
= \tau(m_1)e^{O\left(w/ \tau(m_1)^{1/2}\right)}.
\ee
Since $w= o (\tau(m_1)^{1/2})$, it follows  in particular that
\bel{Sigmaconc}
\tau(m)= \tau(m_1)+ O\left(w\tau(m_1)^{1/2}\right)\sim  \tau(m_1)\quad\mbox{for $m_0\le m\le m_2$}.
\ee
Hence, using~\eqn{lderiv}, we deduce that  for $k=0$ and 2,
\bel{ellell}
|\ell(m_k)-\ell(m_1)| \ge|m_k-m_1|\tau(m_1)/ 3n =K
\ee
where $K= w\tau(m_1)^{1/2}/3$.  (Here we used that $f(b)-f(a)>C(b-a)$ whenever $f'>C$ on $[a,b]$.)

Define the event
$$ A_1= \{L(m_0) < s_1 <L(m_2)- K/2\},
$$
and note that $\ell(m_1)=s_1$ by definition of the inverse function. Observe that since $\beta_j(m_1)\to\infty$, it is necessary that $m_1=o(n^2)$, and hence the same goes for $m_0$ and $m_2$, so we can apply
 Lemma~\ref{l:l}   to $m=m_0$ and $m_2$. This, together with~\eqn{Sigmaconc}, shows that the $\tau$ terms are similar in the two cases, and with \eqn{ellell}, now implies
 \begin{align}\label{PrA1}
1-\pr(A_1) &\le \pr(L(m_0)\ge \ell(m_1)) + \pr(L(m_2)\le \ell(m_1)+K/2)\nonumber\\
&\le\pr(L(m_0)\ge\ell(m_0)+K)+\pr(L(m_2)\le \ell(m_2)-K/2)=O(1/w^2).
\end{align}

Next, define the event
$$
 A_2=\{L(\widetilde M(s_1))-s_1<  K/3 \}.
$$

 We have from~\eqn{secondk} that $\pr(\widetilde B(N-s_1) \ge 2K/3-1 )=O(1/K)$, and hence~\eqn{L-to-s} gives
\bel{PrA2}
 1-\pr(A_{2})=\pr(\widetilde B(N-s_1)/2+1/2\ge K/3) =O(1/K)=O(1/w^2),
\ee
 as $K=o(w^2)$ by the definitions of $w$ and $K$.

  Recall that $\widetilde S(m_0)\le L(m_0)$  and thus the first inequality in $A_1$  implies $\widetilde S(m_0)<s_1$ which is equivalent to  $m_0<\widetilde M(s_1)$. The second inequality in $A_1$, combined with $A_2$, implies
$L(\widetilde M(s_1))<L(m_2)$ and hence, by the monotonicity of $L(m)$,  $\widetilde M(s_1)<m_2$. Thus, from~\eqn{PrA1} and~\eqn{PrA2},
\bel{sandwich}
\pr(m_0<\widetilde M(s_1)<m_2)\ge \pr(A_1\cap A_2) =1- O(w^{-2}).
\ee

Now that we know the number of balls dropped when the number of edges is $s_1$ is highly likely to be between $m_0$ and $m_2$, we show concentration of the vertex degrees at $s_1$ via their concentrations at $m_0$ and $m_2$.
Throughout, we use $j$ to denote an integer satisfying $0\le j\le d-1$.   For $k=0$ and 2, define
 $$
A_{j,k} = \bigg\{\bigg| Y_{\le j}(m_k)- \sum_{i=0}^j\beta_i(m_k) \bigg| < w \sqrt{\beta_j(m_k)}\bigg\}
$$
and notice that for this event  to fail,
$$
\overline{Y_i}(m_k):=|Y_i(m_k) -\beta_i(m_k)|\ge   w\sqrt{\beta_j(m_k)}/(j+1)
$$
must hold for at least one $i$ in the range $0\le i\le j$.

By the assumption $m_1\ge n$ and the definitions of $m_0$ and $w$,  we have $m_0> n/2$ (for $n$ sufficiently large). Thus, recalling that $\mu=m/n$, ~\eqn{tildeYasy} implies
\bel{betasum}
\beta_i(m) = \beta_j(m)[j]_{j-i}/\mu^{j-i}\le (2d)^d\beta_j(m).
\ee
for all $m\ge m_0$ (though all that really matters here is that this is $O(\beta_j(m))$).

 Hence, from Lemma~\ref{l:easy}  with $\alpha= (2d)^{-d/2}w/(j+1)$,  we have for each of $i=0,\ldots, j$, and  $k=0,2$, (note the switch from $\beta_j$ to $\beta_i$)
\bel{Ybeta}
\pr\left( \overline{Y_i}(m_k)\ge \frac{w\sqrt{\beta_j(m_k)}}{j+1}\right)\le\pr\left(\overline{Y_i}(m_k)\ge   \frac{w\sqrt{\beta_i(m_k)}}{(j+1) (2d)^{d/2}} \right)=O(1/w^2),
\ee
  which implies that
$$
\pr(A_{j,k })=1-O( w^{-2})
 $$
 for all $0\le j\le d-1$ and $k=0,2$.

 We next show the following.  (When $j=0$, interpret  $A_{j -1,k }$ as having probability 1.)

 \begin{claim}\lab{Aclaim}  For all $0\le j\le d-1$ and $k=0,2$,
assuming $A_{j,k }$ and  $A_{{j-1},k }$ both hold,   we have
$$
|Y_{j}(m)-  \beta_j(m_1)  | = O( w\beta_j(m_1)  ^{1/2} )
$$
uniformly for all $m$ in the range $m_0\le m\le m_2$.
\end{claim}

\noindent
{\bf Proof\ }
From  the first equation in~\eqn{betaderiv} and~\eqn{betasum} we have the bound $\beta'_i(m)=O(\beta_i(m)/n)$ for all $m\ge m_0$.
Using this, the argument culminating in~\eqn{Sigmaconc} applies also to show that for $k=0$ and $2$
 $$
 \beta_i(m_k) =   \beta_i(m_1) + O\bigg(  w\beta_j(m_1) /\sqrt{\tau(m_1)}\bigg)
  =  \beta_i(m_1) + O(w\beta_i(m_1)  ^{1/2})
$$
where the second step uses $\beta_i\le \tau$.
Thus $A_{j,k }$ implies the same event with each $\beta_i(m_k)$ replaced by  $\beta_i(m_1)$, and $w$ replaced by $O(w)$; that is
\bel{Y-to-sum}
\bigg| Y_{\le j}(m_k)- \sum_{i=0}^j\beta_i(m_1) \bigg| =O(w\beta_j(m_1)  ^{1/2} )
\ee
 for $k=0$ and $ 2$.
Monotonicity of  $Y_{\le j}$ now  yields that~\eqn{Y-to-sum} holds uniformly
for  all $m$ in the range $m_0\le m\le m_2$, which for $j=0$ is the desired claim.
For $j\ge1$,  the claim  follows from the triangle inequality:
 $$|Y_{j}(m)-  \beta_j(m_1)  |\le\bigg| Y_{\le j}(m)- \sum_{i=0}^j\beta_i(m_1) \bigg|+\bigg| Y_{\le j-1}(m)- \sum_{i=0}^{j-1}\beta_i(m_1) \bigg|,$$
recalling, for the error term, that $\beta_{j-1}=O(\beta_j)$. \qed

\medskip

Using our earlier derived bound on $\pr(A_{j,k })$, we deduce from Claim~\ref{Aclaim} that there is a constant $C_0$ such that for $0\le j\le d-1$,
$$
\pr\big(|Y_{j}( m)-  \beta_j(m_1)  |\le  C_0w  \beta_j(m_1)^{1/2}  \mbox{ for   $m_0\le m\le m_2$}\big) =  1-O(w^{-2}).
$$

Also note that  \eqref{Cl22t} can be rewritten as
$|Y_j(m)-D_j(\widetilde S(m))|\le \widetilde B( N-\widetilde S(m))+1$, for all $m$.
 Applying   these two conclusions to $m=\widetilde M(s_1)$, in which case $\widetilde S(m)=s_1$, along with the triangle inequality, gives
$$
\pr\big(|D_{j}( s_1)-  \beta_j(m_1)  |\le  C_0w  \beta_j(m_1)^{1/2}+ \widetilde BN-s_1)+1  \big) = 1-O(w^{-2}).
$$
From~\eqn{secondk} we have  $\pr(\widetilde B(N-s_1)>  w \beta_j(m_1 )^{1/2})=O(1/(w  \beta_j(m_1 )^{1/2} ) = O(w^{-2}) $, since $w=O(\beta_j(m_1 )^{1/2})$.
Thus
$$
\pr\big(|D_{j}( s_1)-  \beta_j(m_1)  |\le  (C_0+2)w  \beta_j(m_1)^{1/2}   \big) = 1-O(w^{-2}).
$$
Replacing $w$ by $(C_0+2)w$ in this gives \eqref{Dconc} in Case 1.

\medskip

\noindent
{\bf Case 2: $m_1 < n$.\/}
\smallskip

 This case is a little trickier than Case 1. The same argument will basically be used with
$L(m)$,  defined in \eqn{Ldef}, rewritten as follows, using the fact that $m=\sum_i iY_i$:
\bel{Ldef-alt}
L(m) = \frac{m}{2}-  \frac12\sum_{i=d+1}^{m}(i-d) Y_i(m).
\ee
Instead of analysing
  $\sum_{i\le j} Y_i$, we will now consider $\sum_{i\ge j} Y_i$.  Since this summation involves an unbounded number of values of $i$ simultaneously, we need  suitable bounds on the deviation  of $Y_i$ from its expected value for arbitrary $i$.  (Lemma~\ref{l:easy}   was adapted to the case that $i$ is fixed.)

We need the following replacement of   Lemma~\ref{l:l}.
 \begin{lemma}\lab{l:l2}
For $1< w =o(\sqrt{\beta_{j}(m_1)}) $ and any $m\le 2n$,
$$
\pr\big(|L(m)-\ell(m)|>w \sqrt{\beta_{d+1} (m)}+1\big) =O(1/w^2).
$$
\end{lemma}
\proof  Luckily, in the summation in~\eqn{Ldef-alt},  values of $i$ greater than  $\log m $ will be neglectable, which can be deduced as follows. From~\eqn{EY} we have $\ex Y_i(m )\le n\mu^i/i! \le   n2^i/i!$. If $m > n^{1/4}$
  then $i>\log m $ implies  $i>\frac14 \log n$, and so $\ex i Y_i(m )=o(1/n^2)$. On the other hand, if  $m \le  n^{1/4}$, then  $\ex i Y_i(m )=O(n^{1-3i/4})$. Since $i\ge d+1\ge 3$, we may now deduce that  $\sum_{i=\lceil \log m\rceil }^{m} \ex iY_i(m)    =o(1/n)$. Similar analysis gives $\sum_{i=\lceil \log m\rceil }^{m} i \beta(m)= o(1/n)$. Hence, by Markov's inequality,
 \bel{big-i}
 \pr\Big(   \sum_{i=\lceil \log m\rceil }^{m} iY_i(m) + i\beta_i(m)  < 1\Big) =1-o(1/n).
\ee

For analysing~\eqn{Ldef-alt}, we are left with estimating  $Y_i(m)$ for $d+1\le i<\log m$. Here we can re-run the calculations that followed~\eqn{EY2} to see that
\bel{EYcase2}
\ex Y_i(m)=\beta_i(m)\left(1+O\left(\frac{i^2}{m}+\frac{i}n+\frac {m}{n^2}\right)\right)
=\beta_i(m)\left(1+O\left(\frac{  i^2  }{m}\right)\right) =O(\beta_i(m))
\ee
  as $m=O(n)$ and $i=O(\log m)$, and similarly
$$
\ex [Y_i(m)]_2=(\ex Y_i(m))^2 \left(1+O\left(\frac{  i^2 }{m}\right)\right)= \ex Y_i(m))^2+O(\ex Y_i(m ))
$$
since  $\ex Y_i(m )\le m\mu^{i-1}/i!\le  m2^{i-1}/i!=O\big(m/ i^2 )\big)$.   It follows (c.f.\ the proof of~\eqn{VarY}) that
 $${\bf Var\,} Y_i(m) =O(\beta_i ).
 $$
 Thus, as in the proof of Lemma~\ref{l:easy},  we have for any $\alpha>0$,  uniformly over $d+1\le i  <\log m$,
\bel{Ybeta2}
\pr( |Y_i(m) -\beta_i(m)|\ge w \sqrt{\beta_{d+1}}/i^3 )=O\left(\frac{i^6 \beta_i}{\beta_{d+1}w^2}\right) = O(w^{-2}/i^2)
\ee
since (recalling  that $m\le 2n$) $\beta_i/\beta_{d+1} \le(2/i)^{i-d-1} =O(i^{-8})$.

  As $\sum_{i\ge 0}\beta_i=n$ and $\sum_{i\ge 0}i\beta_i=m$, we can re-express $\ell$ as
 $$\ell(m) =  \frac{m}{2}-  \frac12\sum_{i=d+1}^{m}(i-d) \beta_i(m).
  $$
Thus, from~\eqn{Ldef-alt},
\bel{Ldiff}
L(m)-\ell(m)  =    \frac12\sum_{i=d+1}^{m}(i-d)( \beta_i(m)-Y_i(m)).
\ee
By~\eqn{big-i}, the terms in the summation in~\eqn{Ldiff} with $i\ge \log m$ sum to  at most $1$ with probability $1-O(w^{-2})$. On the other hand, for $d+1\le i<\log m$,~\eqn{Ybeta2} shows that
$$
\pr\big(|(i-d)( \beta_i(m)-Y_i(m))|\le w \sqrt{\beta_{d+1}(m)}/i^3 \big) =1- O(w ^{-2}/i^2).
$$
 If all of these events hold, then the right-hand side of~\eqn{Ldiff} is at most $w  \sqrt{\beta_{d+1}  (m)}+1 $.
 The sum of the failure probabilities is $O(w^{-2})$. Lemma~\ref{l:l2} follows.
 \qed

  Recalling that $m_1=\ell^{-1}(s_1)$, re-define  $m_0$ and $m_2$ in this case as $ m_0' =m_1 - K' $ and $m_2'=m_1 +K'$ respectively, where $ K' = w\sqrt{\beta_{j}(m_1)}$.

We will consider the event
$$ A'_1= \{L(m'_0) < s_1 <L(m'_2)-  K'/2\},
$$ analogous to the event $A_1$ in Case 1.

From the first relation in~\eqn{pre-conc}  we have
$$
\tau(m)= \tau(m_1)+ O\left(w\tau(m_1)\beta_{j}(m)^{1/2}/n\right)\sim  \tau(m_1)\quad\mbox{for $m'_0\le m\le m'_2$}
$$
as in~\eqn{Sigmaconc} for Case~1.  Following the argument there, we again have from~\eqn{lderiv} that  for $k=0,2$
$$
|\ell(m'_k)-\ell(m_1)|=\Omega(|m'_k-m_1|\tau(m_1)/n)
=\Omega(K)
$$
since $\tau(m_1)=\Omega(n)$ (as $m_1\le n$; here $\beta_0=\Omega(n)$). In this way, this time using Lemma~\ref{l:l2}, we again have
\bel{PrA12}
\pr(A'_1) = 1-O(1/w^2).
\ee
The argument for $A_2=\{L(\widetilde M(s_1))-s_1<  K'/3\}$ is the same as  for $A_2$ in Case~1, and we again obtain~\eqn{sandwich}, which says that with probability $1-O(w^{-2})$ the number of edges reaches $s_1$ when $m$ is between $m'_0$ and  $m'_2$.

For $k=0,2$ define
$$
A'_{j,k} = \bigg\{\bigg| Y_{\ge j}(m'_k)- \sum_{i=j}^{m}\beta_i(m'_k) \bigg| < w \sqrt{\beta_j(m'_k)}\bigg\}.
$$

In this case, we easily have $m'_2<2n$ for large $n$ and hence,  place of~\eqn{betasum},
\bel{betasum2}
\sum_{i=j}^{ m }\beta_i = O(\beta_j).
\ee
Summing~\eqn{Ybeta2} over $d+1\le i  <\log m$, and recalling~\eqn{big-i} for all larger $i$ and Lemma~\ref{l:easy} for $j\le i\le d$,  we get
$$
\pr(A'_{j,k}) = O(w^{-2}).
$$

We can also show the analogue of Claim~\ref{Aclaim}, as follows.

\begin{claim}\lab{Aclaim2}
Assuming $A'_{j,k }$ and  $A'_{j+1,k }$ both hold,   where $  w=o(\sqrt{\beta_j(m_1)})$,  we have
$$
|Y_{j}(m)-  \beta_j(m_1)  | = O\left( w\sqrt{\beta_{j}(m_1)} \right)
$$
uniformly for all $m$ in the range $m'_0\le m\le m'_2$.
\end{claim}

The proof only requires a minor adjustment from Case 1. This time, since $\beta_j(m)$ can be much smaller than $\beta_{j-1}(m)$, for $j\ge 1$  we can only  infer  from~\eqn{betaderiv} and~\eqn{betasum2}, using $\beta_{j-1} (m) /\beta_j (m) =O(n/m)$,   that the derivative of $\beta_j(m)$ satisfies
$$
 \beta'_j(m)  = O(\beta_{j-1} (m)/n)= O(\beta_j(m)/m).
$$
Since $|m'_k-m_1| = K$, we have
$$
 \beta_j(m'_k) =   \beta_j(m_1) + O\bigg(  \beta_j(m_1) w  \sqrt{\beta_j(m_1) }/m\bigg)
  =  \beta_j(m_1) + O(w\beta_j(m_1)  ^{1/2}) ,
$$
since easily $\beta_j (m) \le m$ for $j\ge 1$.

The case $j=0$ is simpler since then we get the   derivative bounded by $\beta_0 (m) /n$ and the same result follows easily.
The rest of the proof of Claim~\ref{Aclaim2} follows the earlier claim's proof almost exactly, and the rest of the proof after that is identical.
We again get~\eqn{Dconc}.


\section{From $Y_i(m)$ to the degrees of vertices at small deficit}\lab{sYD}

For analysing the last stages of the  bin $d$-process,  we find it more convenient to view $D_i$ and $Y_i$ as  functions of the deficit $t$ rather  than the number of balls $m$.
To effect this transition, for $t\ge0$, define $M(t)=\min\{m: T(m)=t\}$, setting $M(t)=\infty$ if this set is empty. Then define
$$
\D(t)=\DD(M(t)).
$$

Let $Y_i^*(t)=Y_i(M(t))$  and $Y_{\le i}^*(t)=Y_{\le i}(M(t))$ where $Y_{\le i}(m)=\sum_{j=0}^i Y_i(m)$.  In this section we first show the concentration of $Y^*_i(t)$  and its factorial moments based on concentration of $Y_i(m)$ established in the previous section. After that, we will further  deduce concentration of $\D(t)$.

\subsection{Sharp concentration of $(Y_i^*(t))_k$}

Throughout we fix  $d\ge 2$ (to avoid trivialities) and $0\le i\le d-1$.   Note that  $Y_{\le i}(m)$  is a decreasing function.  Recall  (cf.\ \eqref{Tdef})  that for every $m$ the random variable $T(m)$ equals the deficit  when $m$ balls have been dropped. In particular, $T(m)$ is decreasing.
 Our first task is to estimate $T(m)$. From~\eqn{Tdef}, the fact that the number of good balls is $\sum_{i=1}^{d-1} i  \DD  (m)$, and the type of simple argument as in the proof of Claim~\ref{determ},
it is easy to see that
 \begin{equation}\lab{Tm}
 T(m)=\frac12\sum_{i=0}^{d-1}(d-i)Y_i(m)+O(B).
 \end{equation}

Let $C>\min\{d,3\}+1$.  For this section, we will make certain assumptions and definitions depending on this fixed value of $C$.   Recall   $\beta_i=\beta_i(n,m)$ as defined in~\eqn{EY}  and assume that $\beta_{d-1} \sim2\log^\gamma n$  for some fixed $\min\{d,3\}+1<\gamma<C$, and $m=\Omega(n)$. Then  Claim~\ref{c:beta} implies that $m\sim n\log n$ and, for $0\le i\le d-1$, that  $\beta_i\sim 2[d-1]_{d-1-i}\log^{\gamma-d+1+i}n$. Moreover, $Y_j(m)$ is sharply concentrated around $\beta_j$ with width  $\beta_j^{3/4}$ for $j\le i$   by Lemma \ref{YtoYtilde}(ii).
 Since $\beta_{j-1}=o(\beta_j)$ when $j>0$, this implies
\bel{Yleq}
 \pr(| Y_{\le i}- \beta_i|\ge 2\beta_i^{3/4})= o(\beta_i^{-K})
 \ee
for every $K>0$.
  Similarly, the most significant term in \eqref{Tm} is that with $i=d-1$. Consequently, noting that $\gamma> 4$, we have,  for every $K>0$, that with probability at least $1-\log^{-K}n$
 $$T(m)=\frac12\beta_{d-1}\left(1+O\left( \frac1{\log n}\right)\right)+O(B).$$
 By Lemma \ref{badballs} with $k>C$, together with Markov's inequality, we have
$$
\pr(B>\log^{2}n) =O(\log^{-C}n),
$$
 and thus
 \begin{equation}\lab{Tm2}
 T(m)=\frac12\beta_{d-1}\left(1+O\left( \frac1{\log n}\right)\right)
 \end{equation}
with probability at least $1-O(\log^{-C}n)$.

Having estimated $T(m)$, we are in a position to estimate $M(t)$ and, through this, functions of $t$ such as $Y^*_i(t)$.

Recall that the above observations assume the constraints on $C$ and $\gamma$ stated in the following lemma.
 \begin{lemma}\lab{l:Ystarasy}
If  $C$ and $\gamma$ are constants satisfying $\min\{d,3\}+1<\gamma<C$ and $t$  is an integer such that  $t \sim \log^\gamma n$, then, for all $0\le i\le d-1$, there exists a constant $C'$ such that with probability $1-O(\log^{-C}n)$ we have
$$
(Y^*_{i}(t))_k=\left(\frac{2[d-1]_{d-1-i}t}{\log^{d-1-i}n}\right)^k\left(1\pm C'\frac1{\sqrt{\log n}}\right).
$$

\end{lemma}
\proof  Set $\beta(m):=\beta(n,m)$ for convenience.
 Let $\widehat M(t)$ be the inverse function of the function $\beta_{d-1}(m)/2$ restricted to the interval  $ I:= \{m:\;\frac12n\log n<m<2n\log n\}$.  (For large $n$, the logarithmic derivative of $\beta_{d-1}$ is negative on $I$, so the inverse exists.)

 Note that $\{ \beta_{d-1}(m)/2\,:\, m\in I\}$ contains values below and above $t$, so $\widehat M(t)$ is defined.  Put $\hat m=\widehat M(t)$, so that
$$
t = \beta_{d-1}(\hat m)/2,
$$
and set $m_1=\hat m-\lfloor n/ \log^{\gamma+1/2}n\rfloor$ and  $m_2=\hat m+\lfloor n/\log^{\gamma+1/2}n\rfloor$.  We  say that an event holds  \emph{with quite  high probability}, abbreviated to \emph{wqhp}, if the event holds with probability $1-O(\log^{-C}n)$. We will show that wqhp
 \begin{equation}\lab{YY}
 T(m_1)>t>T(m_2).
 \end{equation}
For this, we bound the derivative of  $\beta_{d-1}(m)=ne^{-\mu}\mu^{d-1}/(d-1)!$, recalling that $\mu=m/n$.
  On the interval $I$, we have
 $$
 \frac{d \beta_{d-1}(m)}{dm}=\frac{ \beta_{d-1}(m)}n\left(\frac{d-1}{m/n}-1\right)\sim-\frac{\beta_{d-1}(m)}n\sim-\frac{\beta_{d-1}(\hat m)}n.
$$
 Hence
\bean
 \beta_{d-1}(m_1)&=&\beta_{d-1}(\hat m)+(m_1-\hat m)\left(-\frac{\beta_{d-1}(\hat m)}n(1+o(1))\right)\\
&=&2t\left(1+\Theta\left(\frac1{\sqrt{\log n}}\right)\right)
\eean
 since $ \tfrac12\beta_{d-1}(\hat m)=t$ and $m_1-\hat m= O(n/t\sqrt{\log n})$.

Combining with the  symmetric argument for $m_2$, we conclude that
\bel{betaasy}
\beta_{d-1}(m_1)=2t\left(1+\Theta\left(\frac1{\sqrt{\log n}}\right)\right)\quad\mbox{and}\quad  \beta_{d-1}( m_2)=2t\left(1-\Theta\left(\frac1{\sqrt{\log n}}\right)\right).
\ee
In particular, $\beta_{d-1}(m_h)\sim 2\log^\gamma n$ for $h=1$ and $2$. Thus,~\eqn{Tm2} says   that wqhp,  $T(m_h)=\tfrac12\beta_{d-1}(m_h)\left(1+O(1/\log n)\right)$ for $h=1$ and $h=2$. Hence, wqhp
  $$T( m_1)=t\left(1+\Theta\left(\frac1{\sqrt{\log n}}\right)\right)\quad\mbox{and}\quad  T( m_2)=t\left(1-\Theta\left(\frac1{\sqrt{\log n}}\right)\right)$$
  and we have \eqref{YY}.
 It  follows immediately by definition that wqhp
 \begin{equation}\lab{mMm}
 m_1<M(t)<m_2.
 \end{equation}

 Recalling that $t=\tfrac12\beta_{d-1}(\hat m)$,~\eqn{betaasy} implies
 $$(\beta_{d-1}(m_h))_k=(\beta_{d-1}(\hat m))_k\left(1+O\left(\frac1{\sqrt{\log n}}\right)\right)$$
and thus, in view of~\eqn{tildeYasy},
$$(\beta_i(m_h))_k=(\beta_i(\hat m))_k\left(1+O\left(\frac1{\sqrt{\log n}}\right)\right).$$
Hence, by monotonicity of $Y_{\le i}(m)$  and~\eqn{Yleq}, for $0\le i\le d-1$ wqhp
$$(Y_{\le i}(m))_k=(\beta_i(\hat m))_k\left(1+O\left(\frac1{\sqrt{\log n}}\right)\right)$$
for $m_1\le m\le m_2$, and in particular when $m=M(t)$ by \eqref{mMm}. Since  $Y^*_{\le i}(t)=Y_{\le i}(M(t))$ this says that wqhp, for $0\le i\le d-1$,

$$(Y^*_{\le i}(t))_k=(\beta_i(\hat m))_k\left(1+O\left(\frac1{\sqrt{\log n}}\right)\right).$$
 Thus, by~\eqn{tildeYasy},  for each $i$ wqhp
 $$
 (Y^*_{\le i}(t))_k=\left(\frac{2[d-1]_{d-1-i}t}{\log^{d-1-i}n}\right)^k\left(1+O\left(\frac1{\sqrt{\log n}}\right)\right).
 $$
 This finishes the proof of the lemma.  \qed

 \subsection{Sharp concentration of degree counts near the end}

 Here we prove concentration of the degree counts in the bin $d$-process,  and asymptotic formulae for their  moments, at times near the end of the process, as required for our work on the probability of saturation. Theorem~\ref{t:main-degrees} does not provide all the necessary requirements since in some cases we will have $\beta_j(m)\not\to \infty$, and in other cases we need more precise asymptotics. For  convenience we restrict to $t\sim \log^\gamma n$ with a decent lower bound on the constant $\gamma$, though a similar argument applies to interpolated values of $t$.

 \begin{theorem}\lab{t:Ddistrib}
Let $\gamma>\min\{d,3\}+1$ and $t\sim\log^\gamma n$.  Then, setting
$$
f_i(d,t,n)=\frac{2[d-1]_{d-1-i}t}{\log^{d-1-i}n},
$$
for each $0\le i\le d-1$, the following hold for all $k\ge1$.

\begin{enumerate}[(i)]
\item For any fixed  $C>0$, there exists $\Phi(n)=o(1)$ such that with probability $1-O(\log^{-C}n)$
 \begin{equation}\lab{D}
(\D(t))_k=f_i(d,t,n)^k(1\pm \Phi(n)).
\end{equation}
\item We have \begin{equation}\lab{ED}
  \ex((\D(t))_k)\sim f_i(d,t,n)^k.
   \end{equation}
\end{enumerate}

\end{theorem}

\proof Abbreviate $\Bf$ to $B$. By Claim~\ref{determ}, $\D(t)=Y^*_i(t)+O(B)$. Thus,
 $$(\D(t))_k=(Y^*(t)+O(B))_k=\begin{cases}(Y^*_i(t))_k(1+o(1))\quad\mbox{if}\quad B=o(Y^*_i(t))\\O(Y_i^*(t)+B)_k\quad\mbox{always}.\end{cases}$$

 Fix $C>\gamma$. (For smaller $C$ the conclusion follows by monotonicity.) We continue the use of the notation wqhp, with respect to this $C$, from the previous subsection. By applying Lemma \ref{badballs} to the  $j$th moment of $B$ with $j>C$, and Markov's inequality, we get  $B=o(\log^2n)$ wqhp.
 By Lemma~\ref{l:Ystarasy}, wqhp  $Y_i^*(t)\sim f_i(d,t,n)$, where  $ f_i(d,t,n)>\log^2n$  since $\gamma>d+1$.   Hence~\eqref{D} holds wqhp, and we have (i).

We turn to (ii).   Observe that, deterministically,  $\D(t)=O(t)$. Fix $C>(d-1)k$.
 Then, setting $f=f_i(d,t,n)$ and $\Ds=^\neg \hspace{-1.3mm}D(t)$ for convenience, and noting that $(\Ds)_k=O(t^k)$ always,
$$
\ex(\Ds)_k=f^k(1+o(1))\pr\big(|(\Ds)_k-f^k|<\Phi(n)\big)+\pr\big(|(\Ds)_k-f^k|>\Phi(n)\big)O(t^k).
$$
where the latter probability is $O(\log^{-C}n)$ by (i). This implies \eqref{ED}, provided $(t/f)^k=o( \log^{C}n)$, i.e.\ $\log^{(d-1-i)k}n=o( \log^{C}n)$. This holds  for all $0\le i\le d -1$ because $C>(d-1)k$.  \qed


\section{Degree counts late in the process\lab{s:survival}}

Having transferred sharp concentration of the random variables $\DD(m)$ to the variables  $\D(t)$ in the last section, when the deficit of the graph process is small, we no longer need to view it as a bin process. So, we abandon the balls-in-bins terminology altogether and instead focus on analysing the graph process as a function of deficit $t$.

Since concentration is a less effective tool when variables have small values, we will increasingly use estimates of moments, combined with   a special case of Markov's inequality, stating that
 for a nonnegative, integer-valued random variable $X_n$ we have $\pr(X_n\ge1)\le \ex X_n$. In particular, whenever $\ex X_n\to0$ as $n\to\infty$,  a.a.s. $X_n=0$. We will refer to this standard argument as  the \emph{first moment principle}.

In the $d$-process a vertex is called \emph{critical}
  if it has degree at most $d-2$. The next lemma supplements Theorem~\ref{t:Ddistrib} (ii) in an approximate way for times closer to the end of the process. It will be useful in particular for bounding the number of critical vertices.  We use $\Ds_{\le i}$ to denote $\sum_{j=0}^i\Ds_j$.
  \begin{lemma}\lab{l:X}
Fix  integers  $0\le j\le i\le d-2$ and $k \ge1 $, and  a positive  constant  $\gamma$.    For every $t=t(n)\to\infty$ and such that $t<\log^\gamma n$, we have
$$
\ex (\Ds_{\le i}(t))_k = O(\log \log n)^{O(1)} (t\log^{-d+1+i} n) ^k
$$
and
$$
\ex \big(\Ds_{\le i}(t) \Ds_{\le j}(t))= O(\log \log n)^{O(1)}t^2\log^{-2d+2+i+j} n.
$$
\end{lemma}
\proof
Without loss of generality we may assume that $\gamma$ satisfies the lower bound in Theorem~\ref{t:Ddistrib}. Part (ii) of that theorem then implies that both claimed bounds  hold at deficit $t_1=\lceil \log^\gamma n \rceil$ (even without the $\log\log n$ factor).
For the first assertion, note that $(\Ds_{\le i})_k$ counts the ordered $k$-sets of vertices each of degree at most $i$. Moreover, a vertex cannot be of degree at most $i$ when the deficit is $t$ unless it had degree at most $i$ when the deficit was $t_1 \ge t$.  With this in mind, we start by examining such sets of vertices at deficit $t_1$ and consider their evolution until deficit $t$. Condition on the graph at deficit $t_1$ and select  any $k$-set of vertices of degree at most $i$. The probability that they still have degree at most $i$ at deficit $t\le t_1$ and receive some number $a\le dk$ of edges is, by Lemma~\ref{l:survival},   $O(1)\big(\log (t_1 /t )\big)^{a}  ( t/{t_1})^k$.
The expected number of such $k$-sets that still have all members of degree at most $i$ at deficit $t$  is thus
$$\ex (\Ds_{\le i})_k(t_1)\times O(1)\big(\log (t_1 /t )\big)^{a}  ( t/{t_1})^k=O(\log \log n)^{O(1)}(t\log^{-d+1+i} n) ^k,$$
which completes the proof of the first assertion.

 For the second assertion, note that $\Ds_{\le i}  \Ds_{\le j}$ is the number of ordered pairs of  vertices $u$ and $v$ with $d(u)\le i$ and $d(v)\le j$, where possibly $u=v$. This can be expressed as the number of pairs of \emph{distinct} vertices with these degree bounds, plus $\Ds_i$ (as $i\ge j$). An argument very similar to the one above now gives the result.
\qed

Call an edge \emph{unsaturated} if both its endpoints are such, and call it \emph{critical}  if one of its endpoints is unsaturated and the other is  critical.  The next lemma bounds the expected  numbers of unsaturated and critical edges toward the end of the process.

Henceforth, we use $\Lambda$ to stand  for some suitable function of $n$ which is $O(\log\log n)^{O(1)}$, possibly different at each occurrence.

 \begin{lemma}\lab{l:edges}
Let  $ \log^Kn> t \to\infty$ where $K>1$ is fixed.   Then, at deficit $t$
\begin{enumerate}
\item[(i)] The expected number of unsaturated edges   is $ \Lambda t/\log^2 n$.
\item[(ii)] The expected number of pairs of unsaturated edges is $\Lambda t^2/\log^4n$.
\item[(iii)] The expected number of critical edges  is $\Lambda t/\log^3n$.
\end{enumerate}
\end{lemma}
\proof For some fixed $\alpha>2$, take $t_0=\log^{\alpha K}n$.

\smallskip
\noindent (i) Deterministically, since the maximum degree of  the graph is $d$, the number of edges between unsaturated vertices at deficit $t_0$ is $O(t_0)$. By Lemma \ref{l:survival} the probability that any of them survives (i.e.\ its ends remain unsaturated) till deficit $t$ is $\Lambda(t/t_0)^2$. So, the expected number of those which survive till deficit $t$ is
$$\Lambda t^2/t_0=\Lambda/\log^{(2\alpha-2)K}n =o(1/\log^2 n)
$$
since $\alpha>2$.
 This takes care of edges already existing at deficit $t_0$.
For edges created after $t_0$, fix $t_1$ such that $t_0<t_1\le t$. The endvertices of an unsaturated edge created at deficit $t_1$ must have been critical before that edge was added.
    By Lemma \ref{l:X}, the expected number of pairs of critical vertices at deficit $t_1$ is  $\ex \Ds_{\le d-2}(t_1)=\Lambda t_1^2/\log^2n$. Since $U(t_1)=\Theta(t_1)$, the probability that a fixed such pair is selected is $O(1/t_1^2)$.  By Lemma \ref{l:survival}, conditional on such a pair being selected, the probability that they remain unsaturated at deficit $t$ is $\Lambda (t/t_1)^2$. Putting all this together, the expected number of edges between two unsaturated vertices at deficit $t$ that were created after deficit $t_0$ is
$$\sum_{t_1=t}^{t_0}\Lambda t_1^2/\log^2n\times O(1/t_1^2)\times \Lambda (t/t_1)^2=\frac{\Lambda}{\log^2n}\sum_{t_1=t}^{t_0}(t/t_1)^2=\frac{\Lambda t}{\log^2n}.$$

\noindent(ii) We classify the pairs of unsaturated edges present at deficit $t$ into two types: type 1 if at least one of the edges was present at deficit $t_0$, and type 2 otherwise.  By the initial argument in part (i) but with $\alpha$ replaced by $\alpha+4$,  we know already that the expected number of unsaturated  edges at deficit $t$ that were already present at deficit $t_0$ is $\Lambda/\log^{( \alpha+4)K}n$. Each pair of type 1 contains one of these, together with some other edge unsaturated at deficit $t$, of which there are $O(t)=O(\log^K n)$ possibilities. Hence, the number of type 1 pairs present at deficit $t$ is $\Lambda/\log^{K}n=\Lambda/\log^{4}n$.
To create a type 2 pair,  we need  a 4-tuple $(v_1,v_2,v_3,v_4)$ of  critical vertices at $t_0$, without the edges $v_1v_2$ or $v_3v_4$ present at this time (so that those two edges may be added by the time the deficit reaches $t$).  By Lemma \ref{l:X}, the expected number of such 4-tuples is $\Lambda t_0^4/\log^4n$. Given such a 4-tuple, if they determine a type 2 pair with the edge $v_1v_2$ added before $v_3v_4$, the following events must all be satisfied for some $t_1$ and $t_2$ with $t_0>t_1>t_2\ge t$:
\begin{itemize}
\item $A$ -- all four vertices are critical at deficit $t_1+1$ and $e_{t_1}=v_1v_2$,
\item $B$ -- $v_1,v_2$ remain unsaturated and $v_3,v_4$ critical at $t_2+1$, and $e_{t_2}=v_3v_4$,
\item $C$ -- all four remain vertices remain unsaturated at deficit $t$.
\end{itemize}
By Lemma \ref{l:survival}, $\pr(A)=\Lambda(t_1/t_0)^4\times (1/t_1^2)=\Lambda(t_1^2/t_0^4)$, and similarly $\pr(B|A)=\Lambda(t_2^2/t_1^4)$, and $\pr(C|A\cap B)=\Lambda(t/t_2)^4$.
So, the probability that the 4-tuple becomes a pair of unsaturated edges, both created later than $t_0$ is, by the chain formula,
$$
\sum_{t<t_1<t_0}\sum_{t\le t_2<t_1}\pr(A)\pr(B|A)\pr(C|A\cap B)
=\Lambda(t^4/t_0^4)\sum_{t<t_1<t_0}\sum_{t\le t_2<t_1} \frac{1}{t_1^2t_2^2}
=\Lambda(t^2/t_0^4).
$$
Therefore, the expected  number of type 2 pairs edges between unsaturated vertices at $t$ is
$$\Lambda t_0^4/\log^4n\times  t^2/t_0^4 =\Lambda t^2/\log^4n.$$
Combining with the result for type 1 gives (ii).

\smallskip
\noindent (iii) The expected number of critical edges at $t_0$ is  deterministically $O(t_0)$ and, by Lemma \ref{l:survival}, the probability that any one of them remains critical until deficit $t$ is $\Lambda (t/t_0)^2$. So their expected number is $\Lambda t^2/t_0= \Lambda\log^{-3}n$ provided that we choose $\alpha>5$. For the new edges we apply an argument similar to (but simpler than) the argument in (ii). The expected number of pairs of vertices at $t_0$, one of which is critical and the other of degree at most $d-3$, is $\ex \Ds_{ \le d-2}\Ds_{\le d-3}=\Lambda t_0^2/\log^3n$ by Lemma~\ref{l:X}. The event  that a given such non-adjacent pair will end up being a critical edge at deficit $t$ implies that this edge is created at some deficit $t\le t_1<t_0$,    both vertices are unsaturated at deficit $t_1+1$, and  neither vertex becomes saturated later (after $t_1$, up until $t$). Thus, the probability of that event can be bounded, using  Lemma~\ref{l:survival} twice  with $k=2$, by
$$
\sum_{t\le t_1<t_0}\Lambda(t_1/t_0)^2\times 1/t_1^2\times (t/t_1)^2=\Lambda(t/t_0)^2\sum_{t\le t_1<t_0} t_1^{-2}=\Lambda t/t_0^2.
$$
Multiplying by the above bound on $\ex \Ds_{\le d-2}\Ds_{\le d-3}$ shows that the  expected number  of these edges is $\Lambda t/\log^3n$.
\qed

  An interesting consequence of Lemmas \ref{l:edges}(i)  and \ref{l:X}  is  obtained by  the first moment principle.
\begin{corollary}\lab{indept}
For any fixed $\eps>0$, there are a.a.s.\ no  unsaturated  edges  in the $d$-process  at any time after the deficit  reaches  $ \lfloor \log^{1-\eps} n \rfloor$.
\end{corollary}
\proof By Lemma \ref{l:edges}(i) and the first moment principle, there are a.a.s.\ no unsaturated edges at deficit $t =    \lfloor \log^{1-\eps} n \rfloor$. Similarly, using Lemma \ref{l:X}  with $k=1$, at deficit $t$, a.a.s.\ all vertices have degree at least $d-1$. So no new unsaturated  edge can be created afterwards.  \qed

 Next we provide a proof of Theorem \ref{sharpL7}, the most sophisticated argument of the paper. We restate the theorem here with $t$ replaced by $t_1$ just so that $t$ is released to serve as a running variable in the proof.

\medskip

\newtheorem*{thm:main}{Theorem \ref{sharpL7}}
\begin{thm:main}
\textit{For every $j=0,...,d-2$, $k\ge1$, constant $\gamma>0$, and $t_1=t_1(n)<\log^\gamma n$ with $t_1\to\infty$,
 $$\ex\big((\Ds_{j}(t))_k\big)\sim f_{j}(d,t_1,n)^k.
$$}
\end{thm:main}

 \noindent{\bf Proof of Theorem \ref{sharpL7}\ }   We just prove  the case $k=1$ here; the proof of the extension to higher moments is essentially the same, as we explain at the end.

 Let $t_0$ be an integer satisfying $t_0 \sim \log^{d+1+\gamma}n$.
Noting that the degree can only increase as $t$ decreases, there are  just  two contributions to $\ex \Ds_j(t_1)$: vertices of degree $j$  already  at deficit $t_0$, and those of lower degree at $t_0$.
Let
$$
X_1=|\{v:\; \deg_{t_0}(v)=\deg_{t_1}(v)=j\}\quad\mbox{and}\quad X_2=|\{v:\; \deg_{t_0}(v)<j\quad\mbox{and}\quad \deg_{t_1}(v)=j\}
$$
be their numbers respectively. We consider the latter first and show that it is, in fact, negligible.

  By Lemma~\ref{l:X}, $\ex \Ds_{\le  j-1}(t_0)=\Lambda t_0/\log^{d-j}n$. Moreover, for any vertex $v$ and nonnegative integer $i<j$, by Lemma~\ref{l:survival},
$$\pr(\deg_{t_1}(v)=j\;|\;\deg_{t_0}(v)=i)=\Lambda t_1/t_0.$$
Hence
\bel{X2}\ex X_2=\ex \Ds_{\le j-1}(t_0)\Lambda t_0/t_1=\Lambda t_1/\log^{d-j}n.
\ee

\medskip

We use a  much more technical argument  to estimate $\ex X_1$, beginning with some preliminary estimates on the concentration of $U(t)$.
Fix $ \kappa >\kappa>3d+1+2\gamma$,    and $t_0\ge t\ge t_1$.
For $n$ sufficiently large (which we may assume)  we have $t/\log^2t>\kappa$  and hence
$$
\pr\left( \Ds_{\le d-2}(t)>\frac t{\log^2t}\right)=\pr\left((\Ds_{\le d-2}(t))_\kappa>\left(\frac t{\log^2t}\right)_\kappa\right).
$$
By Lemma \ref{l:X}, noting that it applies to all $t$ between $t_1$ and $t_0$ with a suitable adjustment to $\gamma$, with $j=1$,  we have $\ex (\Ds_{\le d-2}(t))_\kappa=\Lambda (t/\log n)^\kappa$. So, the above probability can be   bounded  by Markov's inequality, and noting that $\log^{2\kappa} t = \Lambda$, as follows:
$$
\pr\left((\Ds_{\le d-2}(t))_\kappa>\left(\frac t{\log^2t}\right)_\kappa\right)=\Lambda\left(\frac t{\log n}\right)^\kappa\left(\frac{\log^2t}t\right)^\kappa
=\Lambda \log^{-\kappa}n,
$$
since $\left(  t/\log^2t \right)_\kappa\sim\left(  t/\log^2t \right)^\kappa$. Thus
 \begin{equation}\lab{concl}
\pr\left( \Ds_{\le d-2}(t)>  t/\log^2t \right)=\Lambda \log^{-\kappa}n.
\end{equation}

Recall from Section~\ref{s:distrib} that the number of edges added when the process reaches deficit $t$ is  $s=N-t$, and hence
$$2N-2t=dn-\sum_i (d-i)\Ds_i(t).$$ Thus,
$$
2t=2N-dn+ \sum_i(d-i)\Ds_i(t)\le\sum_i(d-i)\Ds_i(t)\le d\Ds_{\le d-2}(t)+\Ds_{d-1}(t)
$$
and so the number of unsaturated vertices is
$$
U(t)=\sum_{i=0}^{d-1}  \Ds_i(t)\ge \Ds_{d-1}(t)\ge 2t-d\Ds_{\le d-2}(t).
$$
Hence, denoting by $C_t$ the event that  $2t+1\ge U(t)\ge 2t-dt/\log^2t$, it follows from \eqref{concl}  (and the left hand side of \eqref{Ubounds})  that
\bel{Ctbound}
\pr(C_t)=1-\Lambda\log^{-\kappa}n.
\ee

 For a vertex $v$,  our main intermediate goal is to estimate the conditional probability  $\pr( \deg_{t_1}(v)=j\mid \deg_{t_0}(v)=j )$.
  To this end, noting that the event of interest implies that the degree of $v$ remains $j$ for all $t$ between $t_0$ and $t_1$, let $A_{t_0}$ denote the event that $\deg_{t_0}(v)=j$  and let $A_t=A_{t_0}\cap\{\deg_t(v)=j\}$, $t=t_0,\dots,t_1$.   We achieve our goal by inductively showing,  for $t_0\ge t \ge t_1$ (with decreasing $t$),
that
\bel{inductive}
\pr(A_{t-1}|A_t)=1-\frac1t\pm \frac{d}{t\log^2t}
\ee
when $n$ is sufficiently large. It will follow from this  that $\pr(A_t\mid A_{t_0})$ is very well approximated by $t/t_0$, including the case $t=t_1$ (cf. \eqref{ratioComp}).

 We have
\begin{align}
\pr(A_{t-1}|A_t)&=\pr(A_{t-1}|A_t\cap C_t)\pr(C_t|A_t)+\pr(A_{t-1}|A_t\cap \ctbar)\pr(\ctbar|A_t)\nonumber\\
&=\pr(A_{t-1}|A_t\cap C_t)+O\big(\pr(\ctbar|A_t)\big). \lab{onestep}
\end{align}
The last-mentioned probability here is the trickiest to bound. For each vertex $w$, define $F_w$ to be the event that $w$ is unsaturated at deficit $t_0$. Note that $A_t\subseteq A_{t_0}\subseteq F_v$. Hence
\bel{condprob}
\pr(\ctbar|A_t)= \frac{\pr(\ctbar\cap A_t)}{\pr(A_t)}\le\frac{\pr(\ctbar\cap F_v)}{\pr(A_t\cap F_v)}=\frac{\pr(\ctbar\mid F_v)}{\pr(A_t\mid F_v)}.
\ee

Shortly we are going to use vertex symmetry in the $d$-process. Observe that $\ex U_{t_0}=\sum_{w=1}^n\pr(F_w)$ and $\ex (U_{t_0}\mid E)=\sum_{w=1}^n\pr(F_w\mid E)$  for any event $E$. Moreover, $\ex \Ds_j(t_0)=\sum_{w=1}^n\pr(\deg_{t_0}(w)=j)$. By symmetry under relabelling of vertices, with $E=\ctbar$, these identities become, for any fixed $v$,
\bel{symm}
\ex U_{t_0}=n\pr(F_v),\quad \ex( U_{t_0}\mid \ctbar)=n\pr(F_v\mid \ctbar),\quad\mbox{and}\quad \ex \Ds_j(t_0)=nA_{t_0}.
\ee
 To bound the numerator in the right hand side of \eqref{condprob}, we invoke~\eqn{Ubounds}, which in view of \eqref{symm}, yields that for every $w$ we have \newline $\pr(F_w)=\ex U_{t_0}/n\ge 2t_0/(dn)$
and, thus,
\bel{forw}
\pr(\ctbar\mid F_w)\le
\frac{dn}{2t_0}\pr(\ctbar\cap F_w).
\ee
Summing~\eqn{forw}  over  the $n$ different vertices $w$, and using $\pr(\ctbar\mid F_v) = \pr(\ctbar\mid F_w)$ again by symmetry, we get
\bean
n\pr(\ctbar\mid F_v)&\le& \frac{dn}{2t_0}\sum_w \pr(\ctbar\cap F_w)\\
&=&  \frac{dn}{2t_0}\pr(\ctbar)\sum_w \pr( F_w\mid\ctbar)\\
&\overset{\eqref{symm}}{=}&  \frac{dn}{2t_0}\pr(\ctbar)\ex (U_{t_0}\mid\ctbar)\\
&\le&  \frac{dn(2t_0+1)}{2t_0}\pr(\ctbar),
\eean
since $U_{t_0}\le 2t_0 +1$ always by~\eqn{Ubounds}.
Thus
\bel{numerator}
\pr(\ctbar\mid F_v)\le\frac{d(2t_0+1)}{2t_0}\pr(\ctbar)\le2d\pr(\ctbar).
\ee
On the other hand, to estimate the denominator of~\eqn{condprob}, recalling that $A_t\subseteq A_{t_0}\subseteq F_v$, we may write
\bean
 \pr(A_t\mid F_v) =  \pr(A_t\mid A_{t_0}) \cdot  \pr(A_{t_0}\mid F_v)
\eean
Moreover, by \eqref{symm} twice more again,
\begin{align*}
 \pr(A_{t_0}\mid F_v)&=\frac{\pr(A_{t_0})}{\pr(F_v)} =\frac{\pr(A_{t_0})}{\ex(U_{t_0})/n} \ge \frac{n\pr(A_{t_0})}{2t_0+1} = \frac{\ex \Ds_{j}(t_0)}{2t_0+1}\\
&=\Omega(\log^{-d+ 1+j}n),
\end{align*}
the last equality by  Theorem \ref{t:Ddistrib}(ii). Thus,
\bel{denominator}
\pr(A_t\mid F_v) =\pr(A_t\mid A_{t_0})\Omega(\log^{-d+ 1+j}n).
\ee
Substituting estimates \eqref{numerator} and \eqref{denominator} into~\eqn{condprob}, we have
$$
\pr(\ctbar|A_t)= \pr(\ctbar) \frac{O(\log^{d-1-j}n)}{\pr(A_t\mid A_{t_0})}.
$$
To estimate the denominator here,  we  apply~\eqn{inductive} inductively for larger $t$.  Observing that for all $t\le s_1<s_2\le t_0$, the identity $A_{s_1}\cap\cdots\cap A_{s_2}=A_{s_1}$ holds, we obtain, by the chain formula, that
 \bea
\pr(A_t|A_{t_0})&=&
\pr(A_t\cap\cdots\cap A_{t_0-1}|A_{t_0})=
\prod_{s=t+1}^{t_0}\pr(A_{s-1}|A_s) \nonumber\\
&=&\prod_{s=t+1}^{t_0}\left(1-\frac1s\pm \frac{d}{s\log^2s}\right)=
\exp\left\{\sum_{s=t+1}^{t_0}-\frac1s\pm \frac{d}{s\log^2s}+o\left(\frac{1}{s^2}\right)\right\}\nonumber\\
&=&\frac t{t_0}\left(1+O\left(\frac1t\right)+O\left(\frac{ 1 }{\log t}\right)\right)=\frac t{t_0}\left(1+O\left(\frac{ 1 }{\log t}\right)\right).\lab{ratioComp}
\eea
With this estimate of the denominator, recalling that $t\ge t_1\to\infty$  and $t_0\sim\log^{d+1+\gamma}n$, we have
$$
\pr(\ctbar|A_t)= \pr(\ctbar) \frac{t_0}{t}O(\log^{d-1-j}n)
 = \pr(\ctbar)  O(\log^{2d+\gamma}n).
$$
Combining this with~\eqn{onestep} and~\eqn{Ctbound}, we obtain
\begin{align}
\pr(A_{t-1}|A_t)&=\pr(A_{t-1}|A_t\cap C_t)+\Lambda\log^{-\kappa+2d+\gamma}n\nonumber\\
&=\pr(A_{t-1}|A_t\cap C_t)+o\left(\frac1{t\log^2t}\right) \lab{almostThere}
\end{align}
since  $\kappa>3d+1+2\gamma$ and  $t\le t_0 \sim \log^{d+1+\gamma}n$,

The probability that $v$ is incident with the edge added when the deficit is $t$, given $U_t$, is
$$
\frac{U(t)-1-(d-j)}{\binom{U(t)}2-O(U_t)}= \frac2{U(t)}+O\left(\frac1{U(t)^2}\right)
$$
which  is  $\frac1t\pm\frac d{2t\log^2t}(1+o(1))$
 Hence
$$
\pr(A_{t-1}|A_t\cap C_t)=1-\frac1t  \frac{d}{2t\log^2t}(1+o(1)).
$$
For $n$ sufficiently large,  this together with~\eqn{almostThere}
    establishes the inductive hypothesis~\eqn{inductive}.
From~\eqn{ratioComp} with $t=t_1$ we now obtain
$$
\pr(A_{t_1}|A_{t_0})\sim\frac{t_1}{t_0}.
$$

This further implies that
$$
\ex (X_1\mid G_{t_0}) \sim  \frac{t_1}{t_0} \Ds_{j}(t_0).
$$
 By Theorem \ref{t:Ddistrib}(ii), $\ex \Ds_{j}(t_0)\sim \frac{2[d-1]_{d- 1-j}t_0}{\log^{d-1-j}n}$, and so

$$
\ex X_1\sim\frac{t_1}{t_0}\ex \Ds_{j}(t_0)
\sim\frac{t_1}{t_0} \cdot\frac{2[d-1]_{d-1-j}t_0}{\log^{d-1-j}n}
=\frac{2[d-1]_{d-1-j}t_1}{\log^{d-1-j}n}.
$$

Combining this observation with \eqref{X2},
$$\ex \Ds_{j}(t_1)=\ex(X_1+X_2)\sim\frac{2[d-1]_{d-1-j}t_1}{\log^{d-1-j}n}.$$
This proves  the  theorem for $k=1$. For $k\ge2$,  the same proof can be adjusted, by  working with $k$-sets of vertices of degree $j$ and starting with~\eqn{ED} at time $t_0$. We leave these easy adjustments to the reader.\qed

Just the case $k=1$ in Theorem~\ref{sharpL7} provides essentially all the information on distribution of degrees that we need for proving the main results of this paper  on the probability of
saturation,  in the next section.  In particular, we have  the following.

\begin{corollary}\lab{c:Dprob}
For every $j=0,...,d-2$,    constant $\gamma$ satisfying  $0<\gamma<d-1-j$,
and $t  \sim  \log^\gamma n$,
$$
\pr(\Ds_{j}(t )=1)  \sim \frac{2[d-1]_{d-1-j} }{\log^{d-1-j-\gamma}n}.
$$
\end{corollary}
\proof For simplicity, we abbreviate $\Ds_{j}(t)$ to $\Ds$. To  estimate  $\pr(\Ds =1)$ we use the observation that
\bel{EDPD}
\ex \Ds -  \pr( \Ds=1) =\sum_{i\ge 2} i \pr( \Ds=i) \le \sum_{i\ge 2}i(i-1)\pr(\Ds=i)=\ex (\Ds)_2.
\ee
Since  $\gamma< d-1-j$,
$$f_j(d,t,n)=O(t/\log^{d-1-j}n)=o(1).$$
Hence,   Theorem~\ref{sharpL7} gives $\ex (\Ds)_2=o(\ex \Ds)$, which, in view of \eqref{EDPD}, implies  $\pr( \Ds=1)\sim \ex \Ds$. The corollary now follows from  Theorem~\ref{sharpL7} with $k=1$.
\qed

 Theorem~\ref{sharpL7} for general $k$   can be used to provide a more complete   picture of degree distribution in the tail end of the $d$-process, as stated in Theorem \ref{t:Ddistrib2} in the introduction.

\medskip

  \noindent{\bf Proof of Theorem \ref{t:Ddistrib2}\ } This is quite standard. Indeed, it follows immediately from the method of moments for a Poisson random variable  (see~\cite{JLR} for example), using  Theorem~\ref{sharpL7} to estimate the factorial moments in the case $\la$ is bounded, and the result for $\la\to\infty$ follows immediately from this since the limiting probability that $D_j(N-t)=\Ds_{j}(t)=i$ tends to 0 as $t\to\infty$. The two cases can be combined by the subsubsequence principle (as presented in~\cite{JLR}).
 \qed

\medskip

  To round off this section, we easily deduce Corollary \ref{c:lastvertex} from Theorem~\ref{t:Ddistrib2}.

\medskip

\noindent{\bf Proof of Corollary \ref{c:lastvertex}\ }
Consider first $\Ds_0(t)$. This is monotonically decreasing, so the vertices of degree 0 have disappeared iff  $\Ds_0(t) =0$. Thus, this case follows directly from the theorem, and we see that
a.a.s.\ the minimum degree is at least 1 when the deficit is significantly smaller than $\log^{d-1} n$, and definitely by the time it reaches $\log^{d-3/2} n$. Conditional upon this, $\Ds_1(t)$  is monotonic from that point onwards, and repeating the same argument gives the result for all $j<d-1$.
 \qed

We note that the $k$th factorial moments of $\Ds_j$ can also be estimated for $k$ tending to infinity using a similar method, and from this asymptotic normality can be deduced when $\ex \Ds_j\to\infty$ (see~\cite{GW}). However we omit this argument here.

\section{Probability of saturation} \lab{s:sat}
  In this section we prove Theorems \ref{satur} and \ref{t:unsat}. We begin  with a straightforward proof of our old result.

\medskip

\noindent{\bf Proof of Theorem \ref{satur}\ }
Fix $t_0=\lfloor\sqrt{\log n}\rfloor$ and recall that by Corollary~\ref{indept} a.a.s.\  there are  no unsaturated edges  at  deficit $t_0$. On the other hand,
from Lemma \ref{l:X} with $i=d-2$ and $K=1$, we have $\ex \Ds_{\le d-2}(t_0)=\Lambda t_0/\log n$, so by the first moment principle at deficit $t_0$  a.a.s.\ all  unsaturated vertices have degree $d-1$. If this event holds, then after deficit  $t_0$ the process proceeds by a mere addition of disjoint edges and ultimately becomes saturated.  \qed

\medskip

\noindent{\bf Proof of Theorem \ref{t:unsat}\ }  Since the claimed probability is
  $O (1/\log n)$ in the even case and $O(1/\log^2n)$ in the odd case, the main part of the proof is to show that some undesirable events  have even smaller probabilities. We   consider the two assertions separately.

\medskip

\noindent{\bf For  $dn$ even:}
\smallskip

 Fix   $0<\eps<1/2$ and set $t_1=\lceil\log^{\eps}n\rceil$.
 Let $E_1$ be the event that
\begin{enumerate}[(i)]
\item    at deficit $t_1$   there are no unsaturated edges,  and
\item  $\Ds_{d-2}(t_1)\le1$ and $\Ds_{\le d-3}(t_1)=0$.
\end{enumerate}

We first show that
\begin{equation}\label{E1}
\pr(E_1)=1-o(1/\log n).
\end{equation}
 Lemma~\ref{l:edges}(i) says that  the expected number of edges between unsaturated vertices at time $t_1$ is $\Lambda t_1/\log^2n$.  Consequently, by the first moment principle, the probability of sub-event (i)   is $1-\Lambda t_1/ \log^2 n=1-o(1/\log n)$.
Moreover, by Lemma~\ref{l:X}, $\ex \Ds_{\le d-3}(t_1)=\Lambda t_1/\log^2n$ and $\ex (\Ds_{d-2})_2=\Lambda t_1^2/\log^2n$. Hence,  similarly,  $\pr(\Ds_{\le d-3}(t_1)\ge1)=\Lambda t_1/\log^2n$ and $\pr(\Ds_{d-2}(t_1)\ge2)=\pr(((\Ds_{d-2})_2\ge1)= \Lambda t_1^2/\log^2n$.
    Altogether, sub-event (ii) has probability $1-O (t_1^2/ \log^2 n)=o(1/\log n)$ for $\eps<1/2$, and we conclude that \eqref{E1} holds.

\medskip

We next show that $\pr(F\mid E_1)\sim(d-1)/\log n$. This suffices to prove the theorem for $dn$ even, as $\pr(E_1)$ is so close to 1.
Note that $E_1$ implies that for the process not to saturate, at deficit $t_1$ there must be a unique critical vertex, of degree $d-2$, (call it $v$) and all other vertices have degree $d-1$. Moreover, the number of unsaturated vertices  is then odd, because at every step the total remaining degree must be even (and $dn$ is even too). If in the remaining process any added edge is incident with $v$, the process will clearly saturate. So, the process does not saturate if and only if the remaining edges simply match up all the unsaturated vertices other than $v$. Note that if this happens,
\begin{equation}\label{U2}
U(t)=2t-1\quad\mbox{for all}\quad t_1\ge t\ge 1.
\end{equation}

   More precisely, let $A$ be the event that $\Ds_{d-2}(t_1)=1$ and $B$   the event that $\Ds_{d-2}(1)=1$. We have just explained that $F\cap E_1=A\cap B\cap E_1$, so
$$\pr(F\mid E_1)=\pr(A\cap B\mid E_1)=\pr(A|E_1)\pr(B\mid A\cap E_1).$$
We first show that $\pr(B\mid A\cap E_1)=1/(2t_1-1)$.
Indeed, by \eqref{U2},
    $$\pr(B\mid A\cap E_1)=\prod_{t=2}^{t_1}\left(1-\frac{U(t)-1}{\binom{U(t)}2}\right)=\prod_{t=2}^{t_1}\left(1-\frac2{U(t)}\right)=\prod_{t=2}^{t_1}\frac{U(t)-2}{U(t)}
    =\prod_{t=2}^{t_1}\frac{2t-3}{2t-1}=\frac1{2t_1-1}.$$
  By Corollary~\ref{c:Dprob},   $\pr(A)\sim\frac{2(d-1)t_1}{\log n}$, and by \eqref{E1}, $\pr(A\mid E_1)=\pr(A)+o(1/\log n)$. Combined, these observations yield that

$$\pr(F)\sim\pr(F\mid E_1)\sim\frac{2(d-1)t_1}{\log n}\times\frac1{2t_1-1}\sim\frac{d-1}{\log n}.$$

\medskip

\noindent{{\bf For $dn$ odd:}
\smallskip

This is rather harder because the probability of non-saturation is much smaller.
Again let $t_1=\lceil\log ^{\eps}n \rceil$, but $0<\eps<1/3$, and
 let $E_2$ be the event that the following four properties all hold:
\begin{enumerate}[(i)]
\item $\Ds_{\le d-4}(t_1)=0$;
\item $\Ds_{d-2}(t_1)\le 2$, $\Ds_{d-3}(t_1)\le 1$, and  $\Ds_{d-2}(t_1)\Ds_{d-3}(t_1)=0$;
\item     there is at most one unsaturated edge at time $t_1$;
\item  there is no critical edge at time $t_1$.
\end{enumerate}

We first argue that, given $E_2$, the only way for the process not to saturate is by having $\Ds_{d-3}(1)=1$, that is, when a vertex of degree $d-3$ will survive  to the very end.  Similarly to  the even case, let $A$ be the event that $\Ds_{d-3}(t_1)=1$ and $B$  the event that $\Ds_{d-3}(1)=1$. Note that $ B\cap E_2\subseteq A$ by (i) and (ii) in $E_2$.

\begin{claim}\label{theonly}
The event $E_2$ implies that non-saturation occurs if only if $\Ds_{d-3}(1)=1$.  Equivalently, $F\cap E_2=B\cap E_2$.
\end{claim}

\proof The inclusion (implication) $F\cap E_2\supseteq B\cap E_2$ is trivial (as $B\subseteq F$), so we only prove the opposite one.
In $\Gf$, the unsaturated vertices must form a $k$-clique for some $k\ge 1$; otherwise another edge could be added to the graph. If   $k\ge3$, then by (iii)  there is at most one edge between these clique vertices at deficit $t_1$, meaning that they all receive at least one more edge to form the clique by the end. Since they are all finally unsaturated, they must have had  degree at most $d-2$ at deficit $t_1$. However, it is easy to see that (i) and (ii) only permit at most two vertices of degree less than $d-1$ at deficit $t_1$. So, $k\ge 3$ is impossible.

 Assume then that $k=2$.  Due to the total degree parity being even but $dn$ being odd, the two unsaturated vertices must have degrees of different parity in $\Gf$.  Note that by (i) and (ii) it is impossible for one vertex to  have degree at most $d-2$, and another degree at most $d-3$, at any time from deficit $t_1$ onwards. Thus, recalling  (i),  the clique vertices' final degrees  must be $d-1$ and $d-2$, and moreover, at deficit $t_1$ one of them must have been $d-2$ (and hence it was critical) and the other either $d-1$ or $d-2$.  It follows using (iv) that the edge between them was added after deficit $t_1$. However, this contradicts the fact that one of them has final degree $d-2$.  We deduce that there is only one unsaturated vertex in $\Gf$, which must have degree $d-3$, considering parity, and (i), as well as the assumption that the process did not saturate. \qed

\medskip

We will now show that

 \begin{equation}\lab{E2}
\pr(E_2)=1-o(1/\log^2 n).
\end{equation}

 By Lemma \ref{l:X},  $\ex \Ds_{\le d-4}=\Lambda t_1/\log^3n$. Thus, by the first moment principle, (i) fails with probability $O(\log^{\eps-3}n)=o(1/\log^2n)$.

Note that property (ii) is equivalent to restricting the values of the pair $(\Ds_{d-2}(t_1),\Ds_{d-3}(t_1)$ to a set of four options, namely $(0,0),(1,0),(2,0),(0,1)$.
Thus, for  (ii) to fail, since $(i)_j\ge1$ if and only if $i\ge j$, it is required that at deficit $t_1$  at least one of the following holds:
 $(\Ds_{d-2})_3\ge1$,  $(\Ds_{d-3})_2\ge1$, or $\Ds_{d-2}\Ds_{d-3}\ge 1$.

By Lemma \ref{l:X}, we get
\begin{itemize}
\item $\ex (\Ds_{d-2})_3=\Lambda t_1^3/\log^3n$,
\item $\ex (\Ds_{d-3})_2=\Lambda t_1^2/\log^4n$,
\item $\ex \Ds_{d-2}\Ds_{d-3}=\Lambda t_1^2/\log^3n$.
\end{itemize}
All of these are $o(1/\log^2n)$, as $t_1=O(\log^{\eps}n)$ and $\eps<1/3$. Thus, again by the first moment principle, (ii) fails with probability $o(1/\log^2n)$.

For (iii), by Lemma \ref{l:edges}(ii), the expected number of ordered pairs of unsaturated edges  at deficit $t_1$ is $\Lambda t_1^2/\log^4n$. For (iv), by Lemma \ref{l:edges}(iii),  the expected number of critical edges  at deficit $t_1$ is $\Lambda t_1/\log^3n$. This proves \eqref{E2}.

\medskip

With~\eqref{E2} established, it is enough to show that
$$
\pr(F\mid E_2) \sim  \frac{(d-1)(d-2)}{ \log^2 n}.
$$

  Recall, that by Claim \ref{theonly}, $F\cap E_2=A\cap B\cap E_2$, so
$$\pr(F\mid E_2)=\pr(A\cap B\mid E_2)=\pr(A|E_2)\pr(B\mid A\cap E_2).$$
By Corollary~\ref{c:Dprob} and~\eqn{E2},
$\pr(A\mid E_2)\sim\frac{2(d-1)(d-2)t_1}{\log^2 n}$.
Below we will show that
\begin{equation}\label{rest}
\pr(B\mid A\cap E_2)\sim1/2t_1.
\end{equation}
 From these things, we conclude that

$$\pr(F)\sim\pr(F\mid E_2)\sim\frac{2(d-1)(d-2)t_1}{\log n}\times\frac1{2t_1}\sim\frac{(d-1)(d-2)}{\log n}.$$

\medskip

 The rest of the argument consists of showing \eqref{rest}.
For this, we condition on the graph $G_{t_1}$ given that  $E_2$ and $A$ both hold, and suppress this conditioning in the notation.  This implies that there is a   unique vertex $v$ of degree $d-3$  and no vertices other than $v$ are critical at deficit $t_1$.

So now assume that the process fails to saturate. It follows that, for every deficit $t\le t_1$  there must be exactly $2t-2$ vertices of degree $d-1$ and one of degree $d-3$.
 It also follows from  $E_2$   that  when the deficit is $t_1$, at most one edge, which we may call $f$ if it exists, joins two unsaturated vertices, each of which has degree $d-1$, and in view of what we just concluded, the same statement holds for all deficits $t\le t_1$. Hence, for  each  $t\le t_1$,  when the deficit is $t$, conditional upon the current graph $G_t$ with $v$ still of degree $d-3$, the probability of the next edge hitting $v$ is precisely
$$
p_t:= \frac{2t- 2}{(2t-1)(2t-2)/2-\delta_t},
$$
where $\delta_t=1$  if $f$ exists and is still unsaturated, and $\delta_t=0$ otherwise.
Note that $p_t=1/t+O(1/t^2)$ for $t\to\infty$.

Let $t_2 =\lceil  t_1^{1/3}\rceil$ and let $C$ be the event that $\Ds_{d-3}(t_2)=1$, that is,  $v$  retains degree $d-3$  until the deficit reaches $t_2$. Then, suppressing the conditioning on $A\cap E_2$ for convenience,
$$
\pr(C)=\prod_{t=t_2 +1}^{t_1}(1-p_t)
=\exp\left\{\sum_{t=t_2 +1}^{t_1}-\frac1t+O(1/t^2)  \right\}
=\frac{t_2}{t_1} \big(1+O(1/t_2)\big)\sim t_1^{-2/3}.
$$
Further, let $H$ be the event that there is no unsaturated edge at deficit $t_2$ (and thus, afterwards).  Given that $G_{t_1}$ does contain $f$, the probability that both ends of $f$ remain unsaturated until the deficit reaches $t_2$ is, by Lemma~\ref{l:survival}, only $(t_2/t_1)^{2-o(1)}=o(1/t_1)$. So, we have $\pr(H)=1-o(1/t_1)$ and, consequently,
 $$\pr(B)=\pr(B|H)\pr(H)+o(1/t_1)=\pr(B|H)+o(1/t_1).$$
Incorporating conditioning on $C$, we have
$$\pr(B|H)=\pr(B\cap C|H)=\pr(C|H)\pr(B|H\cap C).$$
In turn, $\pr(C|H)=\pr(C)+o(1/t_1)\sim t_1^{-2/3}$, as shown above.

It remains to estimate $\pr(B|H\cap C)$. Here we condition on $G_{t_2}$ not containing an unsaturated edge and $v$ still having degree $d-3$ at deficit $t_2$. In particular, for all $t\le t_2$, we have $\delta_t=0$. Thus,
$$\prod_{t=2}^{t_2}(1-p_t)=\prod_{t=2}^{t_2}\frac{2t-3}{2t-1}=\frac1{2t_2-1}\sim\frac1{2t_2}$$
and, finally, recalling the definition of $t_2$,
$$\pr(B)\sim\frac1{2t_2}\times t_1^{-2/3}\sim\frac1{2t_1}.$$
This completes the proof of Theorem \ref{t:unsat}. \qed

\section{Concluding Remarks}

 Here we briefly mention some open problems and further directions of research based on the balls-in-bins model.

As follows from the proof of Theorem \ref{t:unsat}, the asymptotic probability of non-saturation  coincides with that of having exactly one unsaturated vertex at the end of the process (here, for $dn$ odd a unique vertex of degree $d-1$ is doomed to be saturated). So, how likely is it to see more unsaturated vertices at the end?

\begin{problem} Let $F_k$ be the event that there are $k$ unsaturated vertices in $\Gf$. Estimate $\pr(F_k)$.
\end{problem}
This would require estimating events of smaller probability than we have done. Our result for $dn$ odd required increased accuracy to deal with such events, compared with $dn$ even, and presumably a more careful extension of our argument will provide the necessary higher accuracy.

 We can easily adapt the  techniques of this paper to extend our results to  studying non-saturation probability and degree distribution in $s$-uniform hypergraph $d$-processes,
$s\ge3$. Indeed, the analysis for  hypergraphs is even simpler than for graphs, since, as  was proved in
\cite{GRW}, for $s\ge3$ the \emph{relaxed} hypergraph process (when multiple edges are allowed ---
but still no loops) \aas\ contains no multiple hyperedges. Consequently, arguments such as in the present paper are correspondingly easier.

Our techniques  also  extend  to
 the more restrictive case of \emph{linear} hypergraphs, that is, uniform hypergraphs in which no pair of
  vertices is contained in more than one edge. Here, the probability of bad edges is increased, and the difficulty of arguments is presumably similar to what we have encountered here.

Finally, let us mention that the new balls-in-bins model introduced in this paper opens up the possibility to study other properties of  $d$-processes, as well as  $d$-processes with $d\to\infty$ as $n\to\infty$, since some aspects of the model can be  analysed even in this case.

\subsection*{Acknowledgements}  The authors would like to thank Hania Wormald for helping to facilitate the conditions required for research visits. Also, they would like to thank each  other for their enormous patience and persistence which ultimately lead to the
completion of this marathon project.

\end{document}